\documentclass{article}
\usepackage[dvipsnames,svgnames,x11names]{xcolor}
\usepackage{arxiv}
\usepackage{amsmath}
\numberwithin{equation}{section}
\usepackage[utf8]{inputenc} % allow utf-8 input
\usepackage{hyperref}       % hyperlinks
\usepackage{url}            % simple URL typesetting
\usepackage{booktabs}       % professional-quality tables
\usepackage{amsfonts}       % blackboard math symbols
\usepackage{nicefrac}       % compact symbols for 1/2, etc.
\usepackage{microtype}      % microtypography
\usepackage{cleveref}       % smart cross-referencing
\usepackage{lipsum}         % Can be removed after putting your text content
\usepackage{graphicx}
\usepackage{natbib}
\usepackage{doi}
\usepackage{tikz}
\usepackage{amsthm}
\usepackage{esint}
\tikzset{>=stealth}
\usepackage{bm}
\usepackage{yhmath}
\usetikzlibrary{patterns}
\usepackage{amssymb}
\usepackage{threeparttable}

\title{A new complex variable solution on noncircular shallow tunnelling with reasonable far-field displacement}

\date{\today}
\author{
  \href{https://orcid.org/0000-0002-5143-2714}
  {\includegraphics[scale=0.06]{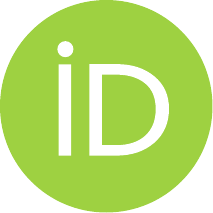}
    \bf {Luo-bin Lin}} \\
  Fujian Provincial Key Laboratory of Advanced Technology and Informatization in Civil Engineering\\
  College of Civil Engineering\\
  Fujian University of Technology\\
  No. 69 Xueyuan Road, Shangjie University Town, Fuzhou, 350118, Fujian, China \\
  \texttt{luobin\_lin@fjut.edu.cn} \\
  \href{https://orcid.org/0000-0002-5583-3734}
  {\includegraphics[scale=0.06]{orcid.pdf}
    \bf {Fu-quan Chen}} \\
  College of Civil Engineering\\
  Fuzhou University\\
  No. 2 Xueyuan Road, Shangjie University Town, Fuzhou, 350108, Fujian, China \\
  \texttt{phdchen@fzu.edu.cn}\\
  \href{}
  {
    \bf {Shang-shun Lin}}\\
  College of Civil Engineering\\
  Fujian University of Technology\\
  No. 69 Xueyuan Road, Shangjie University Town, Fuzhou, 350118, Fujian, China \\
  \texttt{midas2008@126.com} \\
}

\renewcommand{\shorttitle}{}

\hypersetup{
  pdftitle={A new complex variable solution on noncircular shallow tunnelling with reasonable far-field displacement},
  pdfauthor={Luo-bin Lin, Fu-quan Chen, Shang-shun Lin},
  pdfkeywords={Noncircular shallow tunnelling, Far-filed displacement singularity, Stress and displacement, Analytic continuation, Riemann-Hilbert problem},
}

\begin{document}
\maketitle

\begin{abstract}
  A new mechanical model on noncircular shallow tunnelling considering initial stress field is proposed in this paper by constraining far-field ground surface to eliminate displacement singularity at infinity, and the originally unbalanced tunnel excavation problem in existing solutions is turned to an equilibrium one of mixed boundaries. By applying analytic continuation, the mixed boundaries are transformed to a homogenerous Riemann-Hilbert problem, which is subsequently solved via an efficient and accurate iterative method with boundary conditions of static equilibrium, displacement single-valuedness, and traction along tunnel periphery. The Lanczos filtering technique is used in the final stress and displacement solution to reduce the Gibbs phenomena caused by the constrained far-field ground surface for more accurte results. Several numerical cases are conducted to intensively verify the proposed solution by examining boundary conditions and comparing with existing solutions, and all the results are in good agreements. Then more numerical cases are conducted to investigate the stress and deformation distribution along ground surface and tunnel periphery, and several engineering advices are given. Further discussions on the defects of the proposed solution are also conducted for objectivity.
\end{abstract}

% keywords can be removed

\section{Introduction}
\label{sec:intro}

Noncircular tunnels are widely used in civil and mining engineering for the geometrical flexibility of cross section, and noncircular tunnel excavation would generally cause stress and displacement redistribution in nearby geomaterial. For a noncircular deep tunnel, the gravitional gradient can be ignored, and no surface would affect the geometry of geomaterial around the tunnel. Thus, noncircular deep tunnel excavation procedure is generally simplified as a noncircular cavity in an infinite plane, and has been thoroughly investigated~\cite[]{Kargar2014,Lu2014,Lu2015,Wang2019deep_non-circular,fanghuangcheng2021} using the complex variable method~\cite[]{Muskhelishvili1966,Chenziyin_English}.

In contrast, the ground surface limits the geometry of geomaterial to a lower-half plane with an extra surface boundary, and makes the stress and displacement redistribution more complicated even for circular shallow tunnels. The Verruijt conformal mapping~\cite[]{Verruijt1997traction,Verruijt1997displacement} greatly facilitates the application of complex variable method by mapping a lower half plane with a circular cavity onto a unit annulus. With the Verruijt conformal mapping, more complex variable solutions on circular shallow tunnel problems are conducted~\cite[]{Strack2002phdthesis,zhang2018complex,zhang2020complex,kong2019displacement,Ludechun2019,lu2019unified,kong2021analytical}, and these solutions focus on several classic tunnel convergence patterns~\cite[]{Verruijt1996polar}. Many relavent studies using other analytical methods on convergence patterns are also conducted~\cite[]{Sagaseta_image_method,Loganathan_Poulos1998,Bobet2001,Chou2002,park2004,park2005,zhang2017analytical,kong2019displacement,Ludechun2019,lu2019unified,kong2021analytical}.

The convergence patterns would lead to elegant ground surface settlement troughs, but lack neccessary analyses on traction variation along tunnel periphery due to excavation or ground surface. The traction variation along ground surface generally refers to surcharge loads, which are commonly seen in engineering practices, and several complex variable solutions are performed~\cite[]{Wang2018surcharge_twin,Wang2018shallow_surcharge,zhang2021analytical,Self2021APM_in_service,Wanghuaning2021shallow_rigid_lining,zeng2022analytical}.

These solutions related to surcharge loads on ground surface could not convey shallow tunnel excavation procedure with gravational gradient, which is generally simplified as removal of a certain part within gravational lower half plane. Strack~\cite[]{Strack2002phdthesis} introduces additional logarithmic items with singular point in the upper half plane into the complex potentials to cope with the unbalanced resultant imposed along tunnel periphery due to excavation. The modified complex potentials successfully solve the stress of the unbalanced problem, and are further used in many relevant unbalanced problems of circular shallow tunnel excavation~\cite[]{Strack_Verruijt2002buoyancy,Verruijt_Strack2008buoyancy,fang15:_compl,Lu2016,Lu2019new_solution,Self2021APM_shallow_twin}. Furthermore, a new conformal mapping for lower half plane containing a noncircular cavity is proposed by Zeng et al.~\cite[]{Zengguisen2019}, and is soon applied to the unbalanced problems of noncircular shallow tunnel excavations with gravational gradient~\cite[]{Zengguisen2019,lu2021complex} and the noncircular shallow tunnel with surcharge loads~\cite[]{zeng2022analytical}. 

Though the additional logarithmic items of the modified complex potentials successfully lead to correct stress in geomaterial in all the solutions mentioned above, these additional logarithmic items would cause inaccurate displacement instead~\cite[]{Self2020JEM}. The inaccurate displacement consists of two components: one is the displacement singularity at infinity, and the other is the undertermined rigid-body displacement dependent on the artificially selected zero displacement point. Most studies using the modified complex potentials overlook these two inaccurate displacement components or utilize a `translation process` to obtain displacement near the shallow tunnel. To eliminate the displacement singularity at infinity and rigid-body displacement simultaneously, Lin et al.~\cite[]{lin2023reasonable} propose a new mechanical model on shallow tunnel excavation with gravational gradient by tranforming the originally unbalanced problem to an equilbirium one with mixed boundaries along ground surface, and the analytic continuation is used to turn the mixed boundaries into a homogenerous Riemann-Hilbert problem. However, no noncircular tunnels are considered in the solution proposed by Lin et al.~\cite[]{lin2023reasonable}, which limits further possible application of the new mechanical model.

In this paper, a new complex solution is proposed to deal with noncircular shallow tunnel excavation problem considering gravational gradient to simultaneously obtain correct stress and reasonable displacement in geomaterial by eliminating the displacement singularity at infinity.

\section{Problem formation}
\label{sec:problem}

\subsection{Mechanical model of noncircular shallow tunnelling in gravational geomaterial}
\label{sec:problem-1}

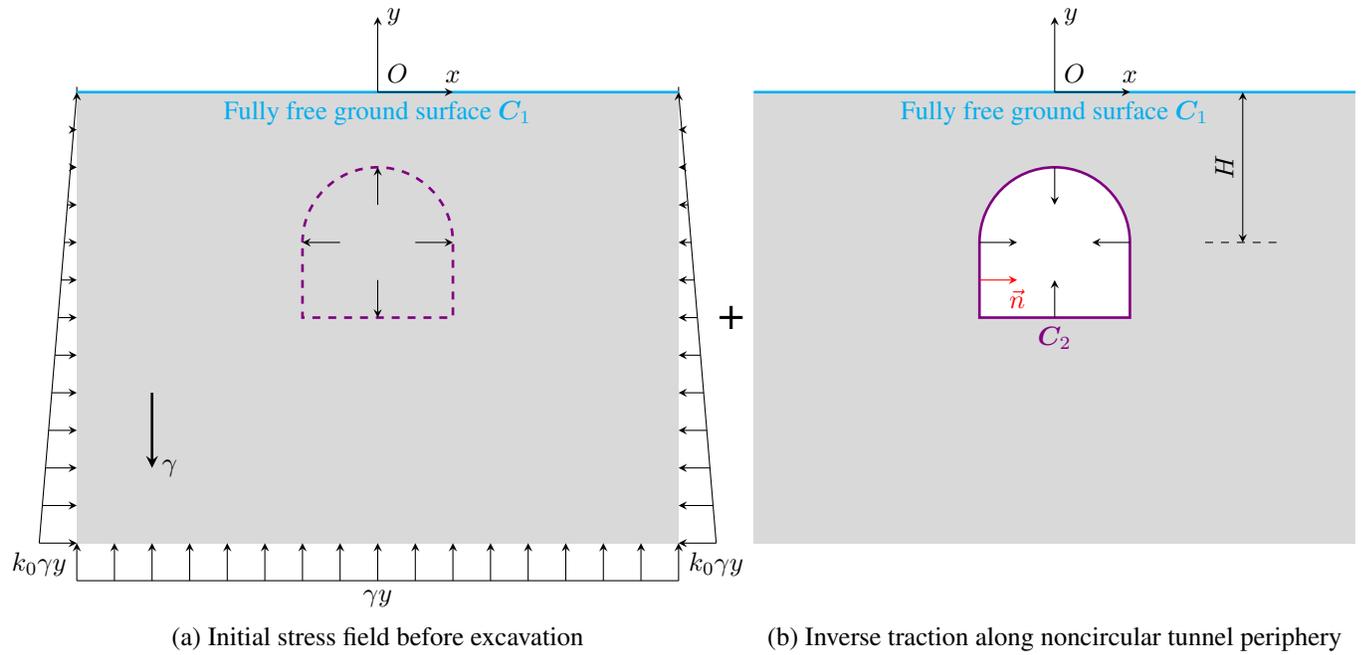
\begin{figure}[htb]
  \centering
  \begin{tikzpicture}
    \fill [gray!30] (0,0) rectangle (8,-6);
    \foreach \x in {0,1,2,...,16} \draw [->] ({\x*0.5},-6.5) -- ({\x*0.5},-6);
    \draw (0,-6.5) -- (8,-6.5);
    \node at (4,-6.5) [below] {$ \gamma y $};
    \foreach \x in {0,1,2,...,12} \draw [->] ({-\x*0.5/12},{-0.5*\x}) -- (0,{-0.5*\x});
    \draw (0,0) -- (-0.5,-6) node [below] {$ k_{0} \gamma y $};
    \foreach \x in {0,1,2,...,12} \draw [->] ({8+\x*0.5/12},{-0.5*\x}) -- (8,{-0.5*\x});
    \draw (8,0) -- (8.5,-6) node [below] {$ k_{0} \gamma y $};
    \draw [cyan, line width = 1pt] (0,0) -- (8,0);
    \draw [->, line width = 1pt] (1,-4) -- (1,-5) node [right] {$ \gamma $};
    \node at (4,0) [cyan, below] {Fully free ground surface $ {\bm C}_{1} $};
    \draw [->] (4,0) -- (5,0) node [above] {$ x $};
    \draw [->] (4,0) -- (4,1) node [right] {$ y $};
    \node at (4,0) [above right] {$ O $};
    \draw [violet, dashed, line width = 1pt] (3,-2) arc [start angle = 180, end angle = 0, radius = 1] -- (5,-3) -- (3,-3) -- (3,-2);
    \draw [->] (4,-1.5) -- (4,-1);
    \draw [->] (4,-2.5) -- (4,-3);
    \draw [->] (3.5,-2) -- (3,-2);
    \draw [->] (4.5,-2) -- (5,-2);
    \node at (4,-7) [below] {(a) Initial stress field before excavation};

    \fill [gray!30] (9,0) rectangle (17,-6);
    \draw [cyan, line width = 1pt] (9,0) -- (17,0);
    \node at (13,0) [cyan, below] {Fully free ground surface $ {\bm C}_{1} $};
    \draw [->] (13,0) -- (14,0) node [above] {$ x $};
    \draw [->] (13,0) -- (13,1) node [right] {$ y $};
    \node at (13,0) [above right] {$ O $};
    \fill [white] (12,-2) arc [start angle = 180, end angle = 0, radius = 1] -- (14,-3) -- (12,-3) -- (12,-2);
    \draw [violet,  line width = 1pt] (12,-2) arc [start angle = 180, end angle = 0, radius = 1] -- (14,-3) -- (12,-3) -- (12,-2);
    \node at (13,-3) [violet, below] {$ {\bm C}_{2} $};
    \draw [->] (13,-1) -- (13,-1.5);
    \draw [->] (13,-3) -- (13,-2.5);
    \draw [->] (12,-2) -- (12.5,-2);
    \draw [->] (14,-2) -- (13.5,-2);
    \draw [->, red] (12,-2.5) -- (12.5,-2.5) node [below] {$ \vec{n} $};
    \draw [dashed] (15,-2) -- (16,-2);
    \draw [<->] (15.5,0) -- (15.5,-2);
    \node at (15.5,-1) [rotate = 90, above] {$ H $};
    \node at (13,-7) [below] {(b) Inverse traction along noncircular tunnel periphery};

    \node at (8.7,-3) {\LARGE{+}};
  \end{tikzpicture}
  \caption{Initial stress field and inverse traction for noncircular shallow tunnel excavation}
  \label{fig:1}
\end{figure}

Zeng et al. \cite[]{Zengguisen2019} propose a mechanical model for noncircular shallow tunnelling in a gravational geomaterial in a lower complex half plane $ z(z = x + {\rm i}y) $, as shown in Fig.~\ref{fig:1}, where the geomaterial $ {\bm \varOmega} $ is homogenerous, isotropic, linearly elastic, of small deformation, and is subjected to the following uniform initial stress field:
\begin{equation}
  \label{eq:2.1}
  \left\{
    \begin{aligned}
      \sigma_{x}^{0} = & \; k_{0} \gamma y \\
      \sigma_{y}^{0} = & \; \gamma y \\
      \tau_{xy}^{0} = & \; 0
    \end{aligned}
  \right.
  , \quad y \leq 0
\end{equation}
where $ \sigma_{x}^{0} $, $ \sigma_{y}^{0} $, and $ \tau_{xy}^{0} $ denote horizontal, vertical, and shear stress components of the initial stress field in the rectangular coordinate system $ xOy $, $ k_{0} $ denotes lateral stress coefficient, $ \gamma $ denotes volumetric weight of geomaterial. The ground surface denoted by $ {\bm C}_{1} $ is fully free from any traction.

Then a noncircular shallow tunnel with buried depth of $h$ is excavated in the geomaterial, whose periphery can be denoted by $ {\bm C}_{2} $. Then tractions in Eq.~(\ref{eq:2.2}) are applied to tunnel periphery to cancel the tractions caused by the initial stress field as shown in Fig.~\ref{fig:1}b:
\begin{equation}
  \label{eq:2.2}
  \left\{
    \begin{aligned}
      X_{i}(S) = & \; - \sigma_{x}^{0}(S) \cdot \cos \langle \vec{n}, x \rangle - \tau_{xy}^{0}(S) \cdot \cos \langle \vec{n}, y \rangle \\
      Y_{i}(S) = & \; - \sigma_{y}^{0}(S) \cdot \cos \langle \vec{n}, y \rangle - \tau_{xy}^{0}(S) \cdot \cos \langle \vec{n}, x \rangle \\
    \end{aligned}
  \right.
\end{equation}
where $ X_{i}(S) $ and $ Y_{i}(S) $ denote the horizontal and vertical traction along the noncircular tunnel periphery $ {\bm C}_{2} $ (denoted by $S$), respectively, $\vec{n}$ denotes the outward direction along the tunnel periphery, as shown in Fig.~\ref{fig:1}b, $ \langle \vec{n},x \rangle $ denotes angle between $ \vec{n} $ and $x$ axis, and $ \langle \vec{n},y \rangle $ denotes angle between $\vec{n}$ and $y$ axis. The tunnel periphery would be traction-free with the sum of the initial stress field in Eq.~(\ref{eq:2.1}) and the tractions in Eq.~(\ref{eq:2.2}). The clockwise integrations of the tractions along tunnel periphery yield the unbalanced resultants along tunnel periphery $ {\bm C}_{2} $:
\begin{equation}
  \label{eq:2.3}
  \left\{
    \begin{aligned}
      F_{x} = & \; \varointclockwise_{{\bm C}_{2}} X_{i}(S) |{\rm d}S| = 0 \\
      F_{y} = & \; \varointclockwise_{{\bm C}_{2}} Y_{i}(S) |{\rm d}S| = 2\pi W \gamma
    \end{aligned}
  \right.
\end{equation}
where $ F_{x} $ and $ F_{y} $ denote horizontal and vertical components of the resultant acting along boundary ${\bm C}_{2}$, respectively; $2\pi W$ denotes the area of the excavated geomaterial, and the detailed expression can be found in Eq.~(\ref{eq:4.17}). The reason of $ F_{x} = 0 $ is the axisymmetry about $y$ axis of the mechanical model.

\subsection{Displacement singularity caused by unbalanced resultant}
\label{sec:problem-2}

Owing to the nonzero unbalanced resultant in Eq.~(\ref {eq:2.3}), the complex potentials $ \varphi (z) $ and $ \psi (z) $ related to the unbalanced resultant can be constructed according to Ref~\cite[]{Zengguisen2019} as
\begin{subequations}
  \label {eq:2.4}
  \begin{equation}
    \label {eq:2.4a}
    \varphi(z) = -\frac{F_{x}+{\rm i}F_{y}}{2\pi(1+\kappa)} \left[ \kappa \ln (z-\overline{z}_{c}) + \ln(z-z_{c}) \right] + \varphi_{0}(z), \quad z \in {\bm \varOmega}
  \end{equation}
  \begin{equation}
    \label {eq:2.4b}
    \psi(z) = -\overline{\varphi}(z)-z\varphi^{\prime}(z) = \frac{F_{x}-{\rm i}F_{y}}{2\pi(1+\kappa)} \left[ \ln(z-\overline{z}_{c}) + \kappa \ln(z-z_{c}) \right] + \psi_{0}(z), \quad z \in {\bm \varOmega}
  \end{equation}
\end{subequations}
where $ z_{c} $ denotes a chosen point within the noncircular tunnel, $ \overline{z}_{c} $ denotes the complex coordinate of the corresponding conjugate point in the upper half plane, $ \varphi_{0}(z) $ and $ \psi_{0}(z) $ denote the single-valued components of the complex potentials, which are always finite in the geomaterial, and
\begin{equation*}
  \psi_{0}(z) = \frac{F_{x}+{\rm{i}}F_{y}}{2\pi(1+\kappa)} \left( \frac{\kappa z}{z-\overline{z}_{c}} + \frac{z}{z-z_{c}} \right) - z\varphi^{\prime}_{0}(z) - \overline{\varphi}_{0}(z)
\end{equation*}

With Eq.~(\ref{eq:2.3}), Eq.~(\ref{eq:2.4}) can be simplified and rewritten in a similar manner of Ref~\cite[]{lin2023reasonable} as
\begin{subequations}
  \label {eq:2.5}
  \begin{equation}
    \label {eq:2.5a}
    \varphi(z) = -\frac{{\rm i}F_{y}}{2\pi} \ln(z-\overline{z}_{c}) + \frac{-{\rm i} F_{y}}{2\pi(1+\kappa)} \ln \frac{z-z_{c}}{z-\overline{z}_{c}} + \varphi_{0}(z), \quad z \in {\bm \varOmega}
  \end{equation}
  \begin{equation}
    \label {eq:2.5b}
    \psi(z) = -\frac{{\rm i}F_{y}}{2\pi}\ln(z-\overline{z}_{c}) + \frac{-{\rm i} \kappa F_{y}}{2\pi(1+\kappa)} \ln \frac{z-z_{c}}{z-\overline{z}_{c}} + \psi_{0}(z), \quad z \in {\bm \varOmega}
  \end{equation}
\end{subequations}
The first items in Eq.~(\ref{eq:2.5}) indicate that a sole upward resultant $ F_{y} $ is imposed at the boundary of the isolated singularity located at coordinate $(0, \overline{z}_{c})$.

The stress and displacement in geomaterial $ {\bm \varOmega} $ can be expressed using the complex potentials in Eq.~(\ref{eq:2.5}) as
\begin{subequations}
  \label{eq:2.6}
  \begin{equation}
    \label{eq:2.6a}
    \left\{
      \begin{aligned}
        \sigma_{y}(z) & + \sigma_{x}(z) = 2 \left[\displaystyle\frac{{\rm d} \varphi(z)}{{\rm d}z} + \overline{\displaystyle\frac{{\rm d} \varphi(z)}{{\rm d}z}}\right] \\
        \sigma_{y}(z) & - \sigma_{x}(z) + 2{\rm i}\tau_{xy}(z) = 2\left[\overline{z} \displaystyle\frac{{\rm d}^{2} \varphi(z)}{{\rm d}z^{2}} + \displaystyle\frac{{\rm d} \psi(z)}{{\rm d}z}\right] \\
      \end{aligned}
    \right.
    , \quad
    z \in {\bm \varOmega}
  \end{equation}
  \begin{equation}
    \label{eq:2.6b}
    g(z) = 2G \left[ u(z) + {\rm i} v(z) \right] = \kappa \varphi(z) - z \overline{\frac{{\rm d} \varphi(z)}{{\rm d}{z}}} - \overline{\psi(z)}, \quad z \in {\bm \varOmega}
  \end{equation}
\end{subequations}
where $ \sigma_{x}(z) $, $ \sigma_{y} (z) $, and $ \tau_{xy} (z) $ denote horizontal, vertical, and shear stress components, respectively; $ u (z) $ and $ v (z) $ denote horizontal and vertical displacement components, respectively; $ G = \frac{E}{2(1+\nu)} $ denotes shear modulus of geomaterial, $ E $ and $ \nu $ denote elastic modulus and Poisson's ratio of geomaterial, respectively, $ \kappa $ denotes the Kolosov coefficient with $ \kappa = 3 - 4\nu $ for plane strain and $ \kappa = \frac{3-\nu}{1+\nu} $ for plane stress.

Substituting Eq.~\ref{eq:2.5} into Eq.~(\ref{eq:2.6}) indicates that the stress components are always finite, while the first items in Eq.~(\ref{eq:2.5}) would result in a unique displacement singularity at infinity in Eq.~(\ref{eq:2.6b}) in a similar manner of Ref~\cite[]{lin2023reasonable}. Thus, we should turn the unbalanced mechanical model above into an equilibrium one by constraining far-field displacement along ground field, similar to Ref~\cite[]{lin2023reasonable}.

\subsection{Boundary conditions}
\label{sec:problem-3}

\begin{figure}[hbt]
  \centering
  \begin{tikzpicture}
    \fill [gray!30] (-5,3) rectangle (3,-3);
    \draw [cyan, line width = 1pt] (-4.5,3) -- (2.5,3);
    \draw [line width = 1.5pt] (-5,3) -- (-4.5,3);
    \draw [line width = 1.5pt] (2.5,3) -- (3,3);
    \fill [pattern = north east lines] (-5,3) rectangle (-4.5,3.2); 
    \fill [pattern = north east lines] (2.5,3) rectangle (3,3.2);
    \fill [white] (-2.5,0) arc [start angle = 180, end angle = 0, radius = 1.5] -- (0.5,-1.5) -- (-2.5,-1.5) -- (-2.5,0);
    \draw [violet, line width = 1pt] (-2.5,0) arc [start angle = 180, end angle = 0, radius = 1.5] -- (0.5,-1.5) -- (-2.5,-1.5) -- (-2.5,0);
    \fill [red] (-4.5,3) circle [radius = 0.05];
    \fill [red] (2.5,3) circle [radius = 0.05];
    \node at (-4.5,3) [below, red] {$ T_{1} $};
    \node at (2.5,3) [below, red] {$ T_{2} $};
    \node at (-2,3) [above, cyan] {$ {\bm C}_{12} $};
    \draw [->] (-1,3) -- (0,3) node [above] {$ x $};
    \draw [->] (-1,3) -- (-1,4) node [right] {$ y $};
    \node at (-1,3) [below] {$ O $};
    \node at (-4.5,3.2) [above] {$ {\bm C}_{11} $};
    \node at (2.5,3.2) [above] {$ {\bm C}_{11} $};
    \node at (-1,-1.5) [above, violet] {$ {\bm C}_{2} $};
    \draw [<->] (1.5,0) -- (1.5,3);
    \draw [dashed] (0.5,0) -- (1.7,0);
    \node at (1.5,1.5) [rotate = 90, above] {$ H $};
    \node at (-1,-4) [below] {(a) Lower half plane containing a noncircular shallow tunnel $ {\bm \varOmega} $};

    \fill [gray!30] (9-0.5,0) circle [radius = 3];
    \draw [cyan, line width = 1pt] (9-0.5,0) circle [radius = 3];
    \fill [pattern = north east lines] (12-0.5,0) -- (12.2-0.5,0) arc [start angle = 0, end angle = 15, radius = 3.2] -- ({9+3*cos(15)-0.5},{3*sin(15)}) arc [start angle = 15, end angle = 0, radius = 3];
    \fill [pattern = north east lines] (12-0.5,0) -- (12.2-0.5,0) arc [start angle = 0, end angle = -15, radius = 3.2] -- ({9+3*cos(-15)-0.5},{3*sin(-15)}) arc [start angle = -15, end angle = 0, radius = 3];
    \draw [line width = 1.5pt] (12-0.5,0) arc [start angle = 0, end angle = 15, radius = 3];
    \draw [line width = 1.5pt] (12-0.5,0) arc [start angle = 0, end angle = -15, radius = 3];
    \fill [white] (9-0.5,0) circle [radius = 1.5];
    \draw [violet, line width = 1pt] (9-0.5,0) circle [radius = 1.5];
    \draw [->] (9-0.5,0) -- (10-0.5,0) node [above] {$ \rho $};
    \draw [->] (9.3-0.5,0) arc [start angle = 0, end angle = 360, radius = 0.3];
    \node at (9.3-0.5,0) [above right] {$\theta$};
    \fill (9-0.5,0) circle [radius = 0.05];
    \node at (9-0.5,0) [below] {$o$};
    \fill [red] ({9+3*cos(15)-0.5},{3*sin(15)}) circle [radius = 0.05];
    \fill [red] ({9+3*cos(-15)-0.5},{3*sin(-15)}) circle [radius = 0.05];
    \node at ({9+3*cos(15)-0.5},{3*sin(15)}) [left, red] {$ t_{2} $};
    \node at ({9+3*cos(-15)-0.5},{3*sin(-15)}) [left, red] {$ t_{1} $};
    \node at (6-0.5,0) [right, cyan] {$ {\bm c}_{12} $};
    \node at (12-0.5,0) [left] {$ {\bm c}_{11} $};
    \node at (9-0.5,-1.5) [above, violet] {$ {\bm c}_{2} $};
    \draw [->] (9-0.5,0) -- ({9-0.5+1.5*cos(-135)},{1.5*sin(-135)}) node [above] {$ r $};
    \node at (9-0.5,-4) [below] {(b) Unit annulus $ {\bm \omega} $};

    \draw [line width = 2pt, ->] (5.2,0) -- (3.2,0);
    \node at (4.2,0) [above] {$ z = z(\zeta) $};
    \node at ({9+3*cos(-135)-0.5},{3*sin(-135)}) [below left] {$ {\bm \omega}^{-} (\rho > 1) $};
    \node at (9-0.5,-2) [below] {$ {\bm \omega}^{+} (\rho < 1) $};
  \end{tikzpicture}
  \caption{Schematic diagram of mixed boundary conditions and conformal mappings of lower half plane containing a noncircular shallow tunnel $ {\bm \varOmega} $ and unit annulus $ {\bm \omega} $}
  \label{fig:2}
\end{figure}
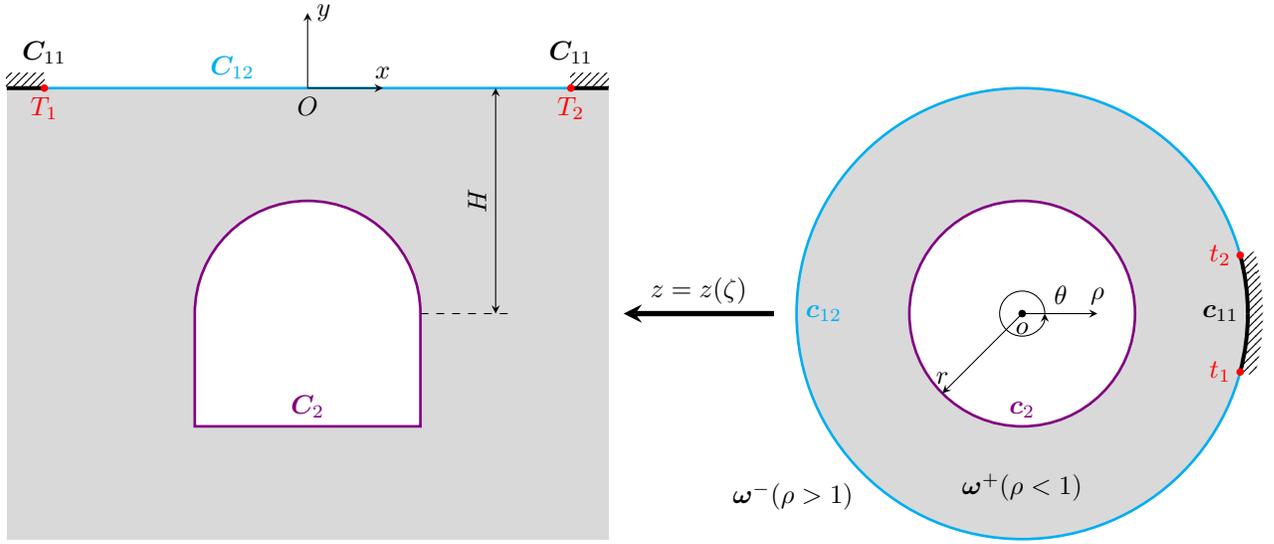

The equilibrium mechanical model is illustrated in Fig~\ref{fig:2}a, where the infinite ground surface is separated into an infinite displacement-constrained segment $ {\bm C}_{11} $ and a finite free segment $ {\bm C}_{12} $ above the noncircular tunnel. Both segments are axisymmetrical, and the intersection points between $ {\bm C}_{11} $ and $ {\bm C}_{12} $ are denoted by $ T_{1} $ and $ T_{2} $, respectively. The following mixed boundary conditions can be constructed in the complex variable manner as
\begin{subequations}
  \label{eq:2.7}
  \begin{equation}
    \label{eq:2.7a}
    u(T) + {\rm i} v(T) = 0, \quad T \in {\bm C}_{11}
  \end{equation}
  \begin{equation}
    \label{eq:2.7b}
    X_{o}(T) + {\rm i} Y_{o}(T) = 0, \quad T \in {\bm C}_{12}
  \end{equation}
\end{subequations}
where $ u(T) $ and $ v(T) $ denote the horizontal and vertical displacement components along boundary $ {\bm C}_{11} $, respectively; $ X_{o}(T) $ and $ Y_{o}(T) $ denote horizontal and vertical tractions along boundary $ {\bm C}_{12} $, respectively. The boundary conditions in Eq.~(\ref{eq:2.7}) are mixed ones.

Eq.~(\ref{eq:2.2}) can be rewritten as
\begin{equation}
  \label{eq:2.2'}
  \tag{2.2'}
  X_{i}(S) + {\rm i} Y_{i}(S) = - k_{0}\gamma y \frac{{\rm d}y}{{\rm d}S} + {\rm i} \gamma y \frac{{\rm d}x}{{\rm d}S}
\end{equation}
The tranformation of Eq.~(\ref{eq:2.2'}) is due to $\cos\langle\vec{n},x\rangle = \frac{{\rm d}y}{{\rm d}{S}}$ and $\cos\langle\vec{n},y\rangle = - \frac{{\rm d}x}{{\rm d}S}$ for outward normal vector $\vec{n}$ and clockwise positive direction of tunnel periphery. Eqs.~(\ref{eq:2.7}) and~(\ref{eq:2.2'}) form the neccessary mathematical expressions for the mixed boundary value problem.

The conformal mapping proposed by Zeng et al.~\cite[]{Zengguisen2019} is then applied to map the geometerial with a noncircular shallow circular tunnel in a lower half plane onto the unit annulus with inner radius of $ r $ (denoted by ${\bm \omega}$ in Fig.~\ref{fig:2}b) via the following mapping function:
\begin{equation}
  \label{eq:2.8}
  z(\zeta) = - {\rm i} a \frac {1 + \zeta} {1 - \zeta} + {\rm i}\sum\limits_{k=1}^{K} b_{k} (\zeta^{k} - \zeta^{-k})
\end{equation}
With Eq.~(\ref{eq:2.8}), the boundaries $ {\bm C}_{11} $, $ {\bm C}_{12} $, and $ {\bm C}_{2} $ in the lower half plane in Fig.~\ref{fig:2}a are respectively mapped by $ {\bm c}_{11} $, $ {\bm c}_{12} $, and $ {\bm c}_{2} $ in the unit annulus in Fig.~\ref{fig:2}b. Boundaries ${\bm c}_{11}$ and ${\bm c}_{12}$ together can be denoted by ${\bm c}_{1}$. Corresponding to region ${\bm \varOmega}$ in the physical plane $z$, region ${\bm \omega}$ in the mapping plane $\zeta$ is also a closure that contains ${\bm c}_{1}$, ${\bm c}_{2}$, and the unit annulus bounded by both boundaries. The regions $ r \leq |\zeta| < 1 $ and $ 1 < |\zeta| \leq r^{-1} $ are denoted by $ {\bm \omega}^{+} $ and $ {\bm \omega}^{-} $, respectively, and $ {\bm \omega} = {\bm \omega}^{+} \cup {\bm c}_{1} $.

It should be addressed that since the infinity point $ z = \infty $ in the lower half plane is mapped onto point $ \zeta = 1 $ in the unit annulus, the infinite boundary $ {\bm C}_{11} $ is mapped onto a finite arc $ {\bm c}_{11} $ correspondingly. The intersection points $ T_{1} $ and $ T_{2} $ are also mapped onto points $ t_{1} $ and $ t_{2} $, respectively. Since $T_{1}$ and $ T_{2} $ are axisymmetrical, $t_{1}$ and $t_{2}$ would also be axisymmetrical. Assume that the poloar angles of points $t_{1}$ and $t_{2}$ would be $-\theta_{0}$ and $\theta_{0}$, respectively, then the horizontal coordinates of points $T_{1}$ and $T_{2}$ are $-x_{0}$ and $x_{0}$, which can be expanded as
\begin{equation}
  \label{eq:2.9}
  x_{0} = \frac{a \sin\theta_{0}}{1-\cos\theta_{0}} - 2\sum\limits_{k=1}^{K} b_{k}\sin k\theta_{0}
\end{equation}
Eq.~(\ref{eq:2.9}) indicates that the range of the free ground surface is dependent on the angle variable $ \theta_{0} $.

Via the conformal mapping in Eq.~(\ref{eq:2.8}), the mixed boundary conditions in Eqs.~(\ref{eq:2.7}) and~(\ref{eq:2.2'}) would be respectively and equivalently mapped by the following ones along corresponding periphery and periphery segements of the unit annulus as
\begin{subequations}
  \label{eq:2.10}
  \begin{equation}
    \label{eq:2.10a}
    u(t) + {\rm i} v(t) = u(T) + {\rm i} v(T) = 0, \quad t = \sigma = {\rm e}^{{\rm i}\theta} \in {\bm c}_{11}
  \end{equation}
  \begin{equation}
    \label{eq:2.10b}
    {\rm e}^{{\rm i}\theta} \frac {z^{\prime}(t)} {|z^{\prime}(t)|} \cdot [\sigma_{\rho}(t) + {\rm i} \tau_{\rho\theta}(t)] = X_{o}(T) + {\rm i} Y_{o}(T) = 0, \quad t = \sigma = {\rm e}^{{\rm i}\theta} \in {\bm c}_{12}
  \end{equation}
  \begin{equation}
    \label{eq:2.10c}
    {\rm e}^{{\rm i}\theta} \frac {z^{\prime}(s)} {|z^{\prime}(s)|} \cdot [\sigma_{\rho}(s) + {\rm i} \tau_{\rho\theta}(s)] = X_{i}(S) + {\rm i} Y_{i}(S) = - k_{0}\gamma y \frac{{\rm d}y}{{\rm d}S} + {\rm i}\gamma y \frac{{\rm d}x}{{\rm d}S}, \quad s = r \cdot \sigma = r {\rm e}^{{\rm i}\theta} \in {\bm c}_{2}
  \end{equation}
\end{subequations}
where points $ t $ and $ s $ in the unit annulus are respectively mapped onto the points $T$ and $S$ in the lower half plane by Eq.~(\ref{eq:2.8}); $ \sigma_{\rho} $ and $ \tau_{\rho\theta} $ denote the radial and shear stress components in the mapping plane, respectively. Now we should solve the mixed boundary value problem in Eq.~(\ref{eq:2.10}).

\section{Problem transformation using analytic continuation}
\label{sec:analyt-cont-riem}

To solve the mixed boundary value problem in Eq.~(\ref{eq:2.10}), the Muskhelishvili's complex potentials $ \varphi(\zeta) $ and $ \psi(\zeta) $ are introduced to express the stress and displacement within geomaterial $ {\bm \omega} $ in the mapping plane $ \zeta = \rho \cdot {\rm e}^{{\rm i}\theta} $ as
\begin{subequations}
  \label{eq:3.1}
  \begin{equation}
    \label{eq:3.1a}
    \sigma_{\theta}(\zeta) + \sigma_{\rho}(\zeta) = 2 \left[ \varPhi(\zeta) + \overline{\varPhi(\zeta)} \right], \quad \zeta \in {\bm \omega} (r \leq \rho \leq 1)
  \end{equation}
  \begin{equation}
    \label{eq:3.1b}
    \sigma_{\rho}(\zeta) + {\rm i}\tau_{\rho\theta}(\zeta) = \varPhi(\zeta) + \overline{\varPhi(\zeta)} - {\rm e}^{-2{\rm i}\theta} \left[ \frac{z(\zeta)}{z^{\prime}(\zeta)}\overline{\varPhi^{\prime}(\zeta)} + \frac{\overline{z^{\prime}(\zeta)}}{z^{\prime}(\zeta)} \overline{\varPsi(\zeta)} \right], \quad \zeta \in {\bm \omega} (r \leq \rho \leq 1)
  \end{equation}
  \begin{equation}
    \label{eq:3.1c}
    g(\zeta) = 2G[u(\zeta)+{\rm i}v(\zeta)] = \kappa \varphi(\zeta) - z(\zeta) \overline{\varPhi(\zeta)} - \overline{\psi(\zeta)}, \quad \zeta \in {\bm \omega} (r \leq \rho \leq 1)
  \end{equation}
\end{subequations}
where
\begin{equation}
  \label{eq:3.2}
  \left\{
    \begin{aligned}
      & \varPhi(\zeta) = \frac{\varphi^{\prime}(\zeta)}{z^{\prime}(\zeta)} \\
      & \varPsi(\zeta) = \frac{\psi^{\prime}(\zeta)}{z^{\prime}(\zeta)} \\
    \end{aligned}
  \right.
\end{equation}

Via analytic contituation, $ \varPhi(\zeta) $ is extended into region $ {\bm \omega} \cup {\bm \omega}^{-} $~\cite[]{Muskhelishvili1966}, and
\begin{equation}
  \label{eq:3.3}
  \psi^{\prime}(\zeta) = \zeta^{-2}\overline{\varphi^{\prime}}\left( \zeta^{-1} \right) + \zeta^{-2}\frac{\overline{z^{\prime}}\left( \zeta^{-1} \right)}{z^{\prime}(\zeta)}\varphi^{\prime}(\zeta) + \frac{\overline{z}\left(\zeta^{-1}\right) z^{\prime\prime}(\zeta)}{[z^{\prime}(\zeta)]^{2}} \varphi^{\prime}(\zeta) - \frac{\overline{z}\left(\zeta^{-1}\right)}{z^{\prime}(\zeta)} \varphi^{\prime\prime}(\zeta), \quad \zeta = \rho \cdot {\rm e}^{{\rm i}\theta} \in {\bm \omega}^{+} (r \leq \rho < 1)
\end{equation}
Consequently, Eq.~(\ref{eq:3.1a}) remains the same, while Eqs.~(\ref{eq:3.1b}) and~(\ref{eq:3.1c}) can be equivalently modified~\cite[]{Muskhelishvili1966,lin2023reasonable} as
\begin{equation}
  \label{eq:3.1b'}
  \tag{3.1b'}
  \begin{aligned}
    \sigma_{\rho}(\zeta) + {\rm i} \tau_{\rho\theta}(\zeta) & = \overline{\zeta}^{2} \left[ \frac{z\left( \overline{\zeta}^{-1} \right)} {z^{\prime}\left( \overline{\zeta}^{-1} \right)} - \frac{1}{\rho^{2}} \frac{z(\zeta)} {z^{\prime}(\zeta)} \right] \overline{\varPhi^{\prime}(\zeta)} + \overline{\zeta}^{2}\left[ \frac{\overline{z^{\prime}(\zeta)}} {z^{\prime}\left( \overline{\zeta}^{-1} \right) } - \frac{1}{\rho^{2}} \frac{\overline{z^{\prime}(\zeta)}} {z^{\prime}(\zeta)} \right] \overline{\varPsi(\zeta)} \\
                                                            & + \varPhi(\zeta) - \varPhi \left( \overline{\zeta}^{-1} \right), \quad \zeta = \rho \cdot \sigma \in {\bm \omega}^{+} (r \leq \rho < 1)
  \end{aligned}
\end{equation}
\begin{equation}
  \label{eq:3.1c'}
  \tag{3.1c'}
  \begin{aligned}
    \frac{{\rm d}g(\zeta)}{{\rm d}\zeta} & =  {\rm e}^{-2{\rm i}\theta} \cdot \left[ z(\zeta) - \rho^{2} \frac{z^{\prime}(\zeta)} {z^{\prime}\left( \overline{\zeta}^{-1} \right)} z\left( \overline{\zeta}^{-1} \right) \right] \overline{\varPhi^{\prime}(\zeta)} + {\rm e}^{-2{\rm i}\theta} \overline{z^{\prime}(\zeta)} \cdot \left[ 1 - \rho^{2} \frac{z^{\prime}(\zeta)} {z^{\prime}\left( \overline{\zeta}^{-1} \right)} \right] \overline{\varPsi(\zeta)} \\
                                                                                                              & + z^{\prime}(\zeta) \left[\kappa  \varPhi(\zeta) + \varPhi\left(\overline{\zeta}^{-1}\right) \right], \quad \zeta = \rho \cdot \sigma \in {\bm \omega}^{+} (r \leq \rho < 1)
  \end{aligned}
\end{equation}
where both of $ \varPhi(\zeta) $ and $ z^{\prime}(\zeta)\varPhi(\zeta) $ should be analytic without any singluarity in the annulus region $ {\bm \omega} \cup {\bm \omega}^{-} $, while $ \varPsi(\zeta) $ and $ z^{\prime}(\zeta)\varPsi(\zeta) $ should be analytic without any singluarity in the unit annulus region $ {\bm \omega} $. The deduction procedure can be found in Ref~\cite[]{lin2023reasonable}, and would not be repeated here. 

When $ \rho \rightarrow 1 $, Eqs.~(\ref{eq:3.1b'}) and~(\ref{eq:3.1c'}) turn to
\begin{subequations}
  \label{eq:3.4}
  \begin{equation}
    \label{eq:3.4a}
    [\sigma_{\rho}(\zeta) + {\rm i} \tau_{\rho\theta}(\zeta)]|_{\rho \rightarrow 1} = \varPhi^{+}(\sigma) - \varPhi^{-}(\sigma)
  \end{equation}
  \begin{equation}
    \label{eq:3.4b}
    \left.\frac{{\rm d}g(\zeta)}{{\rm d}\zeta}\right|_{\rho \rightarrow 1} = z^{\prime}(\sigma) \left[ \kappa \varPhi^{+}(\sigma) + \varPhi^{-}(\sigma) \right]
  \end{equation}
\end{subequations}
where $\varPhi^{+}(\sigma)$ and $\varPhi^{-}(\sigma)$ denote values of $\varPhi(\zeta)$ approaching boundary ${\bm c}_{1}$ from regions ${\bm \omega}^{+}$ and ${\bm \omega}^{-}$, respectively. Comparing Eqs.~(\ref{eq:2.10a}) and~(\ref{eq:2.10b}) with Eqs.~(\ref{eq:3.4b}) and~(\ref{eq:3.4a}) respectively yields
\begin{subequations}
  \label{eq:3.5}
  \begin{equation}
    \label{eq:3.5a}
    \varphi^{\prime+}(\sigma) - \varphi^{\prime-}(\sigma) = 0, \quad \sigma = {\rm e}^{{\rm i}\theta} \in {\bm c}_{12}
  \end{equation}
  \begin{equation}
    \label{eq:3.5b}
    \kappa \varphi^{\prime+}(\sigma) + \varphi^{\prime-}(\sigma) = 0, \quad \sigma = {\rm e}^{{\rm i}\theta}\in {\bm c}_{11}
  \end{equation}
\end{subequations}
where $\varphi^{\prime+}(\sigma)$ and $\varphi^{\prime-}(\sigma)$ denote values of $ \varphi^{\prime}(\zeta) = z^{\prime}(\zeta) \varPhi(\zeta) $ approaching boundary $ {\bm c}_{1} $ from regions $ {\bm \omega}^{+} $ and $ {\bm \omega}^{-} $, respectively. Therefore, the mixed boundary conditions in Eq.~(\ref{eq:2.10}) turn to a homogenerous Riemann-Hilbert problem in Eq.~(\ref{eq:3.5}) with extra constraints in Eq.~(\ref{eq:2.10c}).

\section{Problem solution}
\label{sec:solut-riem-hilb}

\subsection{Laurent series expansion of complex potentials}
\label{sec:solution-riemann-hilbert-1}

Similar to Ref~\cite[]{lin2023reasonable}, the general solution for the Riemann-Hilbert problem in Eq.~(\ref{eq:3.5}) can be expressed as according to Plemelj formula~\cite[]{Muskhelishvili1966}:
\begin{equation}
  \label{eq:4.1}
  \varphi^{\prime}(\zeta) = z^{\prime}(\zeta)\varPhi(\zeta) = X(\zeta) \sum\limits_{n=-\infty}^{\infty} {\rm i} d_{n} \zeta^{n}, \quad \zeta = \rho \cdot \sigma, \quad r \leq \rho \leq r^{-1}
\end{equation}
where
\begin{equation}
  \label{eq:4.2}
  X(\zeta) = (\zeta - {\rm e}^{-{\rm i}\theta_{0}})^{-\frac{1}{2}-{\rm i}\lambda}(\zeta - {\rm e}^{{\rm i}\theta_{0}})^{-\frac{1}{2}+{\rm i}\lambda}, \quad \lambda = \frac{\ln \kappa}{2\pi}
\end{equation}
$ d_{n} $ denote the coefficients to be determined, which should be real owing to the axisymmetry. To solve $ d_{n} $, Eq.~(\ref{eq:4.1}) should be prepared into rational series. $ X(\zeta) $ can be respectivley expanded in regions $ {\bm \omega}^{+} $ and $ {\bm \omega}^{-} $ using Taylor expansions as

\begin{subequations}
  \label{eq:4.3}
  \begin{equation}
    \label{eq:4.3a}
    X(\zeta) = \sum\limits_{k=0}^{\infty} \alpha_{k} \zeta^{k}, \quad \zeta \in {\bm \omega}^{+}
  \end{equation}
  \begin{equation}
    \label{eq:4.3b}
    X(\zeta) = \sum\limits_{k=0}^{\infty} \beta_{k} \zeta^{-k}, \quad \zeta \in {\bm \omega}^{-}
  \end{equation}
\end{subequations}
where 
\begin{equation*}
  \left\{
    \begin{aligned}
      \alpha_{0} = & \; - {\rm e}^{-2\lambda\theta_{0}} \\
      \alpha_{1} = & \; - {\rm e}^{-2\lambda\theta_{0}} \left( \cos\theta_{0} - 2\lambda \sin\theta_{0} \right) \\
      \alpha_{k} = & \; - {\rm e}^{-2\lambda\theta_{0}} \cdot (-1)^{k} \left[ \frac{\mathsf{a}_{k}}{k!} + \frac{\overline{\mathsf{a}}_{k}}{k!} + \sum\limits_{l=1}^{k-1} \frac{\mathsf{a}_{l}}{l!} \cdot \frac{\overline{\mathsf{a}}_{k-l}}{(k-l)!} \right], \quad k \geq 2 \\
      \mathsf{a}_{k} = & \; \prod_{l=1}^{k} \left( \frac{1}{2} - {\rm i}\lambda -l \right) \cdot {\rm e}^{{\rm i}k\theta_{0}}
    \end{aligned}
  \right.
\end{equation*}
\begin{equation*}
  \left\{
    \begin{aligned}
      \beta_{0} = \; & 0 \\
      \beta_{1} = \; & 1 \\
      \beta_{2} = \; & \cos\theta_{0} + 2\lambda \sin\theta_{0} \\
      \beta_{k} = \; & (-1)^{k-1} \left[ \frac{\mathsf{b}_{k-1}}{(k-1)!} + \frac{\overline{\mathsf{b}}_{k-1}}{(k-1)!} + \sum\limits_{l=1}^{k-2} \frac{\mathsf{b}_{l}}{l!} \cdot \frac{\overline{\mathsf{b}}_{k-1-l}}{(k-1-l)!} \right], \quad k \geq 3 \\
      \mathsf{b}_{k} = & \; \prod_{l=1}^{k} \left( \frac{1}{2} - {\rm i}\lambda - l \right) \cdot {\rm e}^{-{\rm i}k\theta_{0}}
    \end{aligned}
  \right.
\end{equation*}

Substituting Eq.~(\ref{eq:4.3}) into Eq.~(\ref{eq:4.1}) yields
\begin{subequations}
  \label{eq:4.4}
  \begin{equation}
    \label{eq:4.4a}
    \varphi^{\prime}(\zeta) = z^{\prime}(\zeta) \varPhi(\zeta) = \sum\limits_{k=-\infty}^{\infty} {\rm i} A_{k} \zeta^{k}, \quad A_{k} = \sum\limits_{n=-\infty}^{k} \alpha_{k-n}d_{n}, \quad \zeta \in {\bm \omega}^{+}
  \end{equation}
  \begin{equation}
    \label{eq:4.4b}
    \varphi^{\prime}(\zeta) = z^{\prime}(\zeta) \varPhi(\zeta) = \sum\limits_{k=-\infty}^{\infty} {\rm i} B_{k} \zeta^{k}, \quad B_{k} = \sum\limits_{n=k}^{\infty} \beta_{n-k}d_{n}, \quad \zeta \in {\bm \omega}^{-}
  \end{equation}
  Since $z^{\prime}(\zeta) \neq 0$ for both regions ${\bm \omega}^{+}$ and ${\bm \omega}^{-}$, Eqs.~(\ref{eq:4.4a}) and~(\ref{eq:4.4b}) can be respectively transformed as
  \begin{equation}
    \label{eq:4.4c}
    \varPhi(\zeta) = \frac{\varphi^{\prime}(\zeta)} {z^{\prime}(\zeta)} = \frac{\sum\limits_{k=-\infty}^{\infty} {\rm{i}}A_{k} \zeta^{k}}{z^{\prime}(\zeta)}, \quad \zeta \in {\bm \omega}^{+}
  \end{equation}
  \begin{equation}
    \label{eq:4.4d}
    \varPhi(\zeta) = \frac{\varphi^{\prime}(\zeta)} {z^{\prime}(\zeta)} = \frac{\sum\limits_{k=-\infty}^{\infty} {\rm{i}}B_{k} \zeta^{k}}{z^{\prime}(\zeta)}, \quad \zeta \in {\bm \omega}^{-}
  \end{equation} 
\end{subequations}
These four above equations analytically indicate that $ \varPhi(\zeta) $ and $ z^{\prime}(\zeta)\varPhi(\zeta) $ are analytic within the annulus region $ {\bm \omega} \cup {\bm \omega}^{-} $, which is identical to the claim in Eqs.~(\ref{eq:3.1b'}) and~(\ref{eq:3.1c'}). Substituting Eqs.~(\ref{eq:4.4b}), (\ref{eq:4.4a'}), and (\ref{eq:2.8}) into Eq.~(\ref{eq:3.3}) with notations of $ \overline{z}(\zeta^{-1}) = z(\zeta) $ and $ \overline{z^{\prime}}(\zeta^{-1}) = -\zeta^{2} z^{\prime}(\zeta) $ yields
\begin{equation}
  \label{eq:4.5}
  \psi^{\prime}(\zeta) = -\sum\limits_{k=-\infty}^{\infty} {\rm i}B_{k}\zeta^{-k-2} + \left\{ \frac{z(\zeta)\cdot z^{\prime\prime}(\zeta)}{[z^{\prime}(\zeta)]^{2}} - 1 \right\} \sum\limits_{k=-\infty}^{\infty} {\rm i}A_{k}\zeta^{k} - \frac{z(\zeta)}{z^{\prime}(\zeta)} \sum\limits_{k=-\infty}^{\infty} {\rm i}kA_{k} \zeta^{k-1}, \quad \zeta \in {\bm {\omega}^{+}}
\end{equation}
Eq.~(\ref{eq:4.5}) can be rewritten into a more compact form according to Eq.~(\ref{eq:3.3}) as
\begin{equation}
  \label{eq:4.5'}
  \tag{4.5'}
  \psi^{\prime}(\zeta) = -\sum\limits_{k=-\infty}^{\infty} {\rm i}B_{-k-2}\zeta^{k} - \left[ \frac{z(\zeta)}{z^{\prime}(\zeta)} \varphi^{\prime}(\zeta) \right]^{\prime}, \quad \zeta \in {\bm {\omega}^{+}}
\end{equation}

\subsection{Static equilbirium and displacement single-valuedness}
\label{sec:solution-riemann-hilbert-2}

Owing to static equilibrium, the absolute value of the constraining force along the constrained segment of the ground surface should be equal to that of the unbalanced resultant caused by excavation along tunnel periphery, but their direction should be opposite. Such mechanical relationship can be expressed by the following expression:
\begin{equation}
  \label{eq:4.6}
  \int_{{\bm C}_{11}} \left[ X_{o}(T)+{\rm i}Y_{o}(T) \right]|{\rm d}T| = - \varointclockwise_{{\bm C}_{2}} \left[ X_{i}(S)+{\rm i}Y_{i}(S) \right]|{\rm d}S|
\end{equation}
where $|{\rm d}T|$ and $|{\rm d}S|$ denote length increments along boundaries ${\bm C}_{11}$ and ${\bm C}_{2}$ in the physical plane, respectively. Since no traction is applied along boundary ${\bm C}_{12}$ in Eq.~(\ref{eq:2.10b}), the left-hand side of Eq.~(\ref{eq:4.6}) can be modified as
\begin{equation}
  \label{eq:4.7}
  \int_{{\bm C}_{11}} \left[ X_{o}(T)+{\rm i}Y_{o}(T) \right]|{\rm d}T| = \int_{{\bm C}_{1}} \left[ X_{o}(T)+{\rm i}Y_{o}(T) \right]|{\rm d}T| = \ointctrclockwise_{{\bm c}_{1}} {\rm e}^{{\rm i}\theta}z^{\prime}({\rm e}^{{\rm i}\theta})[\sigma_{\rho}({\rm e}^{{\rm i}\theta})+{\rm i}\tau_{\rho\theta}({\rm e}^{{\rm i}\theta})]{\rm d}\theta
\end{equation}
where $ |{\rm d}T| = |z^{\prime}({\rm e}^{{\rm i}\theta})|\cdot |{\rm i}{\rm e}^{{\rm i}\theta}| \cdot |{\rm d}\theta| = |z^{\prime}({\rm e}^{{\rm i}\theta})|{\rm d}\theta $ for counter-clockwise length increment in the mapping plane. Substituting Eqs.~(\ref{eq:3.5a}) and (\ref{eq:4.4}) into Eq.~(\ref{eq:4.7}) yields
\begin{equation}
  \label{eq:4.8}
  \ointctrclockwise \left(\sum\limits_{k=-\infty}^{\infty} A_{k}{\rm e}^{{\rm i}k\theta} - \sum\limits_{k=-\infty}^{\infty} B_{k}{\rm e}^{{\rm i}k\theta}\right) \cdot {\rm i}{\rm e}^{{\rm i}\theta} {\rm d}\theta = 2\pi{\rm i}(A_{-1}-B_{-1})
\end{equation}
Based on Eq.~(\ref{eq:2.4}), the right-hand side of Eq.~(\ref{eq:4.6}) can be written as
\begin{equation}
  \label{eq:4.9}
  -\varointclockwise_{{\bm C}_{2}} \left[ X_{i}(S)+{\rm i}Y_{i}(S) \right]|{\rm d}S| = - \varointclockwise_{{\bm C}_{2}} X_{i}(S)|{\rm d}S| - {\rm i}\varointclockwise_{{\bm C}_{2}} Y_{i}(S)|{\rm d}S| = - F_{x} - {\rm i} F_{y} = - 2\pi {\rm i}W \gamma
\end{equation}
Eqs.~(\ref{eq:4.8}) and (\ref{eq:4.9}) should be equal, and we have
\begin{equation}
  \label{eq:4.10}
  A_{-1}-B_{-1} = - W \gamma
\end{equation}

The displacement in geomaterial should be single-valued. Substituting Eqs.~(\ref{eq:4.4a}), (\ref{eq:4.4a'}), and (\ref{eq:4.5'}) into Eq.~(\ref{eq:3.1c}) yields
\begin{equation}
  \label{eq:4.11}
  \begin{aligned}
    g(\zeta) = & \; \kappa \int z^{\prime}(\zeta)\varPhi(\zeta) {\rm d}\zeta - z(\zeta) \overline{\varPhi(\zeta)} - \int \overline{z^{\prime}(\zeta)\varPsi(\zeta)} {\rm d}\overline{\zeta} + {\rm i} C_{0} \\
    = & \; {\rm i}(\kappa A_{-1} + B_{-1}) {\rm Ln}{\rm e}^{{\rm i}\theta} + {\rm i}(\kappa A_{-1} - B_{-1})\ln\rho + {\rm single \; valued \; items}
  \end{aligned}
\end{equation}
where ${\rm Ln}$ denotes the multi-valued natural logarithm sign, $C_{0}$ denotes the integral constant to ensure $g(1)=0$ to satisfy the displacement boundary condition in Eq.~(\ref{eq:2.10a}). To guarantee displacement single-valuedness of Eq.~(\ref{eq:4.11}), the coefficient of the $ {\rm{Ln}} $ item should be zero. Thus, we have
\begin{equation}
  \label{eq:4.12}
  \kappa A_{-1} + B_{-1} = 0
\end{equation}
Finally, Eqs.~(\ref{eq:4.10}) and (\ref{eq:4.12}) yield
\begin{subequations}
  \label{eq:4.13}
  \begin{equation}
    \label{eq:4.13a}
    A_{-1} = \frac{-W\gamma}{1+\kappa}
  \end{equation}
  \begin{equation}
    \label{eq:4.13b}
    B_{-1} = \frac{\kappa W\gamma}{1+\kappa}
  \end{equation}
\end{subequations}

\subsection{Boundary condition expansion of Eq.~(\ref{eq:2.10c})}
\label{sec:solution-riemann-hilbert-3}

To uniquely determine the unknown coefficients $d_{n}$ in Eq.~(\ref{eq:4.1}), the traction boundary condition in Eq.~(\ref{eq:2.10c}) should be also satisfied in the equivalent path integral form. Note that the integral path is clockwise to always keep the geomaterial on the left hand. The first equation in Eq.~(\ref{eq:2.10c}) can be equivalently modified as
\begin{equation}
  \label{eq:4.14}
  - {\rm i} \int [X_{i}(S)+{\rm i}Y_{i}(S)]\cdot|{\rm d}S| = - {\rm i} \int {\rm e}^{{\rm i}\theta} \frac{z^{\prime}(s)}{|z^{\prime}(s)|} [\sigma_{\rho}(s)+{\rm i}\tau_{\rho\theta}(s)] \cdot |z^{\prime}(s)| \cdot |{\rm d}s| = \int z^{\prime}(s)[\sigma_{\rho}(s)+{\rm i}\tau_{\rho\theta}(s)] {\rm d}{s}
\end{equation}
where $ |{\rm d}s| = |r{\rm d}{\rm e}^{{\rm i}\theta}| = r|{\rm d}\theta| = - r{\rm d}\theta $ for clockwise length increment in mapping plane. Simultaneously, the second equilibrium in Eq.~(\ref{eq:2.10c}) can be modified as
\begin{equation}
  \label{eq:4.15}
  - {\rm i} \int [X_{i}(S)+{\rm i}Y_{i}(S)]\cdot|{\rm d}S| = {\rm i} \int \left( - k_{0}\gamma y\frac{{\rm d}y}{{\rm d}S} + {\rm i}\gamma y\frac{{\rm d}x}{{\rm d}S}\right){\rm d}S = - {\rm i} k_{0}\gamma \int y{\rm d}y - \gamma \int y{\rm d}x
\end{equation}
where $|{\rm d}S| = - {\rm d}S$ for clockwise length increment in physical plane.

Both Eqs.~(\ref{eq:4.14}) and (\ref{eq:4.15}) should be prepared into rational series in the mapping plane. Substituting Eqs.~(\ref{eq:4.4}), (\ref{eq:4.5'}) and~(\ref{eq:3.1b}) into Eq.~(\ref{eq:4.14}) yields
\begin{equation}
  \label{eq:4.14'}
  \tag{4.14'}
  \begin{aligned}
    \int z^{\prime}(s)[\sigma_{\rho}(s)+{\rm i}\tau_{\rho\theta}(s)] {\rm d}{s} 
    = & \; \int \varphi^{\prime}(s){\rm d}s + \frac{z(s)}{\overline{z^{\prime}(s)}}\overline{\varphi^{\prime}(s)} + \int \overline{\psi^{\prime}(s)}{\rm d}\overline{s} \\
    = & \; \sum\limits_{k=1}^{\infty} {\rm i} A_{-k-1} \frac{r^{-k}}{-k} \sigma^{-k} + \sum\limits_{k=1}^{\infty} {\rm i}A_{k-1}\frac{r^{k}}{k}\sigma^{k} + \sum\limits_{k=1}^{\infty} {\rm i}B_{k-1} \frac{r^{-k}}{-k} \sigma^{k} + \sum\limits_{k=1}^{\infty} {\rm i}B_{-k-1} \frac{r^{k}}{k}\sigma^{-k} \\
    - & \; \frac{z(r\sigma)-\overline{z(r\sigma)}}{\overline{z^{\prime}(r\sigma)}} \sum\limits_{k=-\infty}^{\infty} {\rm i}A_{-k} r^{-k} \sigma^{k} + {\rm i}(A_{-1}+B_{-1})\ln{r} + {\rm i}(A_{-1}-B_{-1}){\rm Ln}\sigma + {\rm i} C_{a}
  \end{aligned}
\end{equation}
where $C_{a}$ denotes the integral constant to be determined and should be real due to symmetry.

Considering the conformal mapping in Eq.~(\ref{eq:2.6}), the real variables $x$ and $y$ along boundary ${\bm C}_{2}$ in Eq.~(\ref{eq:4.15}) can be written as
\begin{subequations}
  \label{eq:4.16}
  \begin{equation}
    \label{4.16a}
    x = - \frac{{\rm i}a}{2}\left( \frac{1+r\sigma}{1-r\sigma} - \frac{\sigma+r}{\sigma-r} \right) + {\rm i}\sum\limits_{k=1}^{K} \frac{b_{k}}{2} (r^{k}+r^{-k}) (\sigma^{k}-\sigma^{-k})
  \end{equation}
  \begin{equation}
    \label{4.16b}
    y = -\frac{a}{2}\left( \frac{1+r\sigma}{1-r\sigma} + \frac{\sigma+r}{\sigma-r} \right) + \sum\limits_{k=1}^{K} \frac{b_{k}}{2} (r^{k}-r^{-k})(\sigma^{k}+\sigma^{-k})
  \end{equation}
  \begin{equation}
    \label{eq:4.16c}
    {\rm d}x = \left\{ -{\rm i}ar\left[\frac{1}{(1-r\sigma)^{2}}+\frac{1}{(\sigma-r)^{2}}\right] + {\rm i}\sum\limits_{k=1}^{K} \frac{k b_{k}}{2} (r^{k}+r^{-k})(\sigma^{k-1}+\sigma^{-k-1}) \right\}{\rm d}\sigma
  \end{equation}
\end{subequations}
Substituting Eq.~(\ref{eq:4.16}) into Eq.~(\ref{eq:4.15}) yields
\begin{equation}
  \label{eq:4.15'}
  \tag{4.15'}
  - {\rm i} k_{0}\gamma \int_{A}^{B} y{\rm d}y - \gamma \int_{A}^{B} y{\rm d}x = {\rm i}\gamma F(\sigma) + {\rm i} \gamma \sum\limits_{j=1}^{3} \left[E_{j}(\sigma) + D_{j} \cdot {\rm Ln}\sigma \right]
\end{equation}
where
\begin{equation*}
  F(\sigma) = -\frac{k_{0}}{2} \left[ -\frac{a}{2}\left( \frac{1+r\sigma}{1-r\sigma} + \frac{\sigma+r}{\sigma-r} \right) + \sum\limits_{k=1}^{K} \frac{b_{k}}{2} (r^{k}-r^{-k})(\sigma^{k}+\sigma^{-k}) \right]^{2}
\end{equation*}
Detailed expressions of $ D_{j} $ and $ E_{j}(\sigma) $ can be found in Appendix \ref{sec:expansion}.

According to Eqs.~(\ref{eq:4.14}) and (\ref{eq:4.15}), Eqs.~(\ref{eq:4.14'}) and Eq.~(\ref{eq:4.15'}) should be equal. To guarantee such a requirement, the multi-valued components of $ {\rm Ln}{\rm e}^{{\rm i}\theta} $ should be eliminated simultaneously, and we obtain
\begin{equation}
  \label{eq:4.10'}
  \tag{4.10'}
  A_{-1} - B_{-1} = \sum\limits_{j=1}^{3} D_{j}\gamma
\end{equation}
Eqs.~(\ref{eq:4.10}) and~(\ref{eq:4.10'}) should be the same, thus, we have
\begin{equation}
  \label{eq:4.17}
  W = -\sum\limits_{j=1}^{3} D_{j}
\end{equation}
Eq.~(\ref{eq:4.17}) gives the area of the excavated geomaterial in Eq.~(\ref{eq:2.3}). The deduction above further analytically emphasizes the mechanical fact that stress and traction in the geomaterial should be single-valued. The remaining items of Eqs.~(\ref{eq:4.14'}) and~(\ref{eq:4.15'}) give
\begin{equation}
  \label{eq:4.18}
  \begin{aligned}
    & \sum\limits_{k=1}^{\infty} A_{-k-1} \frac{r^{-k}}{-k} \sigma^{-k} + \sum\limits_{k=1}^{\infty} A_{k-1}\frac{r^{k}}{k}\sigma^{k} + \sum\limits_{k=1}^{\infty} B_{k-1} \frac{r^{-k}}{-k} \sigma^{k} + \sum\limits_{k=1}^{\infty} B_{-k-1} \frac{r^{k}}{k}\sigma^{-k} \\
    - & \; \frac{z(r\sigma)-\overline{z(r\sigma)}}{\overline{z^{\prime}(r\sigma)}} \sum\limits_{k=-\infty}^{\infty} A_{-k} r^{-k} \sigma^{k} + (A_{-1}+B_{-1})\ln{r} +  C_{a} = \gamma F(\sigma) + \gamma \sum\limits_{l=1}^{3} E_{l}(\sigma)
  \end{aligned}
\end{equation}

\subsection{Approximate solution}
\label{sec:Approximate-solution}

We should solve $ d_{n} $ in Eq.~(\ref{eq:4.1}) using Eqs.~(\ref{eq:4.13}) and~(\ref{eq:4.18}). Eq.~(\ref{eq:4.18}) can be rewritten as
\begin{equation}
  \label{eq:4.18'}
  \tag{4.18'}
  \begin{aligned}
    & \sum\limits_{k=1}^{\infty} A_{-k-1} \frac{r^{-k}}{-k} \sigma^{-k} + \sum\limits_{k=1}^{\infty} A_{k-1}\frac{r^{k}}{k}\sigma^{k} + \sum\limits_{k=1}^{\infty} B_{k-1} \frac{r^{-k}}{-k} \sigma^{k} + \sum\limits_{k=1}^{\infty} B_{-k-1} \frac{r^{k}}{k}\sigma^{-k} \\
    + & \; \sum\limits_{k=-\infty}^{\infty} I_{k} \sigma^{k} \sum\limits_{k=-\infty}^{\infty} A_{-k} r^{-k} \sigma^{k} + (A_{-1}+B_{-1})\ln{r} +  C_{a} = \sum\limits_{k=-\infty}^{\infty} J_{k} \sigma^{k}
  \end{aligned}
\end{equation}
where
\begin{subequations}
  \label{eq:4.19}
  \begin{equation}
    \label{eq:4.19a}
    \sum\limits_{k=-\infty}^{\infty} I_{k} \sigma^{k} \sim - \frac{z(r\sigma)-\overline{z(r\sigma)}}{\overline{z^{\prime}(r\sigma)}}
  \end{equation}
  \begin{equation}
    \label{eq:4.19b}
    \sum\limits_{k=-\infty}^{\infty} J_{k} \sigma^{k} \sim \gamma F(\sigma) + \gamma \sum\limits_{j=1}^{3} E_{j}(\sigma)
  \end{equation}
\end{subequations}
The symbol $\sim$ denotes that the latter expression can be simulated by the former one, and the computation procedure can be found in Appendix \ref{sec:coefficients}.

Comparing the coefficients of the same order of $ \sigma $ with substitution of Eq.~(\ref{eq:4.4}) yields
\begin{subequations}
  \label{eq:4.20}
  \begin{equation}
    \label{eq:4.20a}
    \sum\limits_{n=k+1}^{\infty}\alpha_{n-k-1}d_{-n} = r^{2k} B_{-k-1} + kr^{2k}\sum\limits_{l=-\infty}^{\infty} r^{l}I_{l}A_{k+l} - kr^{k}J_{-k}, \quad k \geq 1
  \end{equation}
  \begin{equation}
    \label{eq:4.20b}
    \sum\limits_{n=k}^{\infty}\beta_{n-k+1}d_{n} = r^{2k} A_{k-1} + k \sum\limits_{l=-\infty}^{\infty} r^{l}I_{l}A_{-k+l} - kr^{k}J_{k}, \quad k \geq 1
  \end{equation}
  \begin{equation}
    \label{eq:4.20c}
    C_{a} = -\sum\limits_{l=-\infty}^{\infty} r^{l}I_{l}A_{l} - (A_{-1}+B_{-1})\ln{r} + J_{0}
  \end{equation}
\end{subequations}
The sum of Eq.~(\ref{eq:4.20b}) starts from $k$ instead of $(k-1)$, because $ \beta_{0} = 0 $ in Eq.~(\ref{eq:4.3b}). Eqs.~(\ref{eq:4.20a}) and~(\ref{eq:4.20b}) respectively constrain $ d_{-n} (n \geq 2) $ and $ d_{n} (n \geq 1) $, while $ d_{0} $ and $ d_{-1} $ are not constrained. Now Eq.~(\ref{eq:4.13}) should be used and transformed as
\begin{equation}
  \label{eq:4.13a'}
  \tag{4.13a'}
  \sum\limits_{n=1}^{\infty} \alpha_{n-1}d_{-n} = \frac{-W\gamma}{1+\kappa}
\end{equation}
\begin{equation}
  \label{eq:4.13b'}
  \tag{4.13b'}
  \sum\limits_{n=0}^{\infty} \beta_{n+1}d_{n} = \frac{\kappa W\gamma}{1+\kappa}
\end{equation}
Subsequently, $ d_{-n} (n \geq 2) $ can be determined by Eqs.~(\ref{eq:4.20a}) and~(\ref{eq:4.13a'}), while $ d_{n} (n \geq 1) $ can be determined by Eqs.~(\ref{eq:4.20b}) and~(\ref{eq:4.13b'}). Thus, the following approximate solution can be used.

Assume that $ d_{n} $ can be accumulated as
\begin{equation}
  \label{eq:4.21}
  d_{n} = \sum\limits_{q=0}^{\infty} d_{n}^{(q)}
\end{equation}
When $q=0$, the initial values of $ d_{-n}^{(0)} (n \geq 1) $ and $ d_{n} (n \geq 0) $ can be respectively determined by Eqs.~(\ref{eq:4.20a}) and~(\ref{eq:4.13a'}), (\ref{eq:4.20b}) and (\ref{eq:4.13b'}) as
\begin{subequations}
  \label{eq:4.22}
  \begin{equation}
    \label{eq:4.22a}
    \left\{
      \begin{aligned}
        \sum\limits_{n=1}^{\infty} \alpha_{n-1}d_{-n}^{(0)} = & \; \frac{-W\gamma}{1+\kappa} \\
        \sum\limits_{n=k+1}^{\infty}\alpha_{n-k-1}d_{-n}^{(0)} = & \; - kr^{k}J_{-k}, \quad k \geq 1 \\
      \end{aligned}
    \right.
  \end{equation}
  \begin{equation}
    \label{eq:4.22b}
    \left\{
      \begin{aligned}
        \sum\limits_{n=0}^{\infty} \beta_{n+1}d_{n}^{(0)} = & \; \frac{\kappa W\gamma}{1+\kappa} \\
        \sum\limits_{n=k}^{\infty}\beta_{n-k+1}d_{n}^{(0)} = & \; - kr^{k}J_{k}, \quad k \geq 1 \\
      \end{aligned}
    \right.
  \end{equation}
\end{subequations}
For iteration $q \geq 0 $, $A_{k}^{(q)}$ and $B_{k}^{(q)}$ can be computed according to Eq.~(\ref{eq:4.4}) as
\begin{subequations}
  \label{eq:4.23}
  \begin{equation}
    \label{eq:4.23a}
    A_{k}^{(q)} = \sum\limits_{n=-k}^{\infty} \alpha_{n+k}d_{-n}^{(q)}
  \end{equation}
  \begin{equation}
    \label{eq:4.23b}
    B_{k}^{(q)} = \sum\limits_{n=k}^{\infty} \beta_{n-k}d_{n}^{(q)}
  \end{equation}
\end{subequations}
For the next iteration $q+1$, the iterative value of $ d_{-n}^{(q+1)} (n \geq 1) $ and $ d_{n}^{(q+1)} (n \geq 0) $ can be respectively computed via Eqs.~(\ref{eq:4.20a}) and~(\ref{eq:4.20b}) as
\begin{subequations}
  \label{eq:4.24}
  \begin{equation}
    \label{eq:4.24a}
    \left\{
      \begin{aligned}
        \sum\limits_{n=1}^{\infty} \alpha_{n-1}d_{-n}^{(q+1)} = & \; 0 \\
        \sum\limits_{n=k+1}^{\infty}\alpha_{n-k-1}d_{-n}^{(q+1)} = & \; r^{2k} B_{-k-1}^{(q)} + kr^{2k}\sum\limits_{l=-\infty}^{\infty} r^{l}I_{l}A_{k+l}^{(q)}, \quad k \geq 1 \\
      \end{aligned}
    \right.
  \end{equation}
  \begin{equation}
    \label{eq:4.24b}
    \left\{
      \begin{aligned}
        \sum\limits_{n=0}^{\infty} \beta_{n+1}d_{n}^{(q+1)} = & \; 0 \\
        \sum\limits_{n=k}^{\infty}\beta_{n-k+1}d_{n}^{(q+1)} = & \; r^{2k} A_{k-1}^{(q)} + k \sum\limits_{l=-\infty}^{\infty}r^{l}I_{l}A_{-k+l}^{(q)}, \quad k \geq 1 \\
      \end{aligned}
    \right.
  \end{equation}
\end{subequations}
When the following threshold is reached, the iteration may stop to accumulate the final results of $ d_{n} $ via Eq.~(\ref{eq:4.21}):
\begin{equation}
  \label{eq:4.25}
  \max| d_{n}^{(q)} | \leq \varepsilon
\end{equation}
where $\varepsilon$ is a small numeric, and can be chosen as $ \varepsilon = 10^{-16} $, which is the default value of double precision of programming code {\tt{FORTRAN}}. The constant $ C_{a} $ in Eq.~(\ref{eq:4.20c}) can also be solved using the similar iteration procedure above without any technical difficulty, but it would not affect the stress and displacement distribution in the geomaterial, and is consequently not neccessary to solve. 

To pratically obtain numerical solution, Eq.~(\ref{eq:4.1}) should be truncated into $ 2N+1 $ items, and Eq.~(\ref{eq:4.4}) should be truncated as
\begin{equation}
  \label{eq:4.4a'}
  \tag{4.4a'}
  \varphi^{\prime}(\zeta) = \sum\limits_{k=-N}^{N} {\rm i} A_{k} \zeta^{k}, \quad A_{k} = \sum\limits_{n=-k}^{N} \alpha_{n+k}d_{-n}, \quad \zeta \in {\bm \omega}^{+}
\end{equation}
\begin{equation}
  \label{eq:4.4b'}
  \tag{4.4b'}
  \varphi^{\prime}(\zeta) = \sum\limits_{k=-N}^{N} {\rm i} B_{k} \zeta^{k}, \quad B_{k} = \sum\limits_{n=k}^{N} \beta_{n-k}d_{n}, \quad \zeta \in {\bm \omega}^{-}
\end{equation}
Correspondingly, all the rest series in the solution procedure in Eqs.~(\ref{eq:4.5})-(\ref{eq:4.25}) should be truncated as well. Since the truncated series to iteratively obtain $ d_{n}^{(q)} $ have very similar forms to Eqs.~(\ref{eq:4.22})-(\ref{eq:4.24}), these series are not presented here to reduce repeating expressions.

\section{Stress and displacement in geomaterial}
\label{sec:stress-displ-geom}

The solution procedure in last section gives the neccessary coefficients $ d_{n} $ in Eq.~(\ref{eq:4.1}), which would further determine the first derivatives of complex potentials within the geomaterial region $ {\bm{\omega}} $ in Eqs.~(\ref{eq:4.4a}) and~(\ref{eq:4.5}). The stress distribution for arbitrary point $ \zeta = \rho\sigma (r \leq \rho \leq 1) $ in the mapping plane can be obtained according to Eqs.~(\ref{eq:3.1a}) and~(\ref{eq:3.1b}) as
\begin{subequations}
  \label{eq:5.4}
  \begin{equation}
    \label{eq:5.1a}
    \begin{aligned}
      z^{\prime}(\rho\sigma) [\sigma_{\theta}(\rho\sigma) + \sigma_{\rho}(\rho\sigma)]
      = 2 \left[ \sum\limits_{k=-\infty}^{\infty} {\rm{i}}A_{k}\rho^{k}\sigma^{k} - \frac{z^{\prime}(\rho\sigma)} {\overline{z^{\prime}(\rho\sigma)}} \sum\limits_{k=-\infty}^{\infty} {\rm{i}}A_{k}\rho^{k}\sigma^{-k} \right]
    \end{aligned}
  \end{equation}
  \begin{equation}
    \label{eq:5.1b}
    \begin{aligned}
      z^{\prime}(\rho\sigma) [\sigma_{\rho}(\rho\sigma) + {\rm i}\tau_{\rho\theta}(\rho\sigma)] =
      & \; \sum\limits_{k=-\infty}^{\infty} {\rm{i}}A_{k}\rho^{k}\sigma^{k} - \frac{z^{\prime}(\rho\sigma)} {\overline{z^{\prime}(\rho\sigma)}}\sum\limits_{k=-\infty}^{\infty} {\rm{i}}A_{k}\rho^{k}\sigma^{-k} -  \sum\limits_{k=-\infty}^{\infty} {\rm{i}}A_{k}\rho^{k}\sigma^{-k-2} - \sum\limits_{k=-\infty}^{\infty} {\rm{i}}B_{k}\rho^{-k-2}\sigma^{k} \\
      & \; - \frac{z(\rho\sigma)-\overline{z(\rho\sigma)}} {\overline{z^{\prime}(\rho\sigma)}} \cdot \frac{\overline{z^{\prime\prime}(\rho\sigma)}} {\overline{z^{\prime}(\rho\sigma)}} \sum\limits_{k=-\infty}^{\infty} {\rm{i}}A_{k}\rho^{k}\sigma^{-k-2} + \frac{z(\rho\sigma)-\overline{z(\rho\sigma)}} {\overline{z^{\prime}(\rho\sigma)}} \sum\limits_{k=-\infty}^{\infty} {\rm{i}}kA_{k}\rho^{k-1}\sigma^{-k-1} \\
    \end{aligned}
  \end{equation} 
  \begin{equation}
    \label{eq:5.1c}
    \begin{aligned}
      g(\rho\sigma) = 2G[u(\rho\sigma)+{\rm i}v(\rho\sigma)] =
      & \; {\rm{i}}\sum\limits_{k=1}^{\infty}\left( \kappa A_{k-1}\rho^{k} + B_{k-1}\rho^{-k} \right)\frac{\sigma^{k}}{k} - {\rm{i}}\sum\limits_{k=1}^{\infty}\left( \kappa A_{-k-1}\rho^{-k} + B_{-k-1}\rho^{k} \right)\frac{\sigma^{-k}}{k} + {\rm{i}}C_{0} \\
      & \; + \frac{z(\rho\sigma)-\overline{z(\rho\sigma)}} {\overline{z^{\prime}(\rho\sigma)}} \cdot {\rm{i}} \sum\limits_{k=-\infty}^{\infty} A_{k}\rho^{k}\sigma^{-k} + {\rm{i}}(\kappa A_{-1} - B_{-1})\ln{\rho}
    \end{aligned}
  \end{equation}
\end{subequations}

Eq.~(\ref{eq:5.1b}) can be transformed as
\begin{equation}
  \label{eq:5.2}
  z^{\prime}(\rho\sigma) [\sigma_{\rho}(\rho\sigma) + {\rm i}\tau_{\rho\theta}(\rho\sigma)] = \frac{1}{\rho} \frac{{\rm{d}} \left\{ \int z^{\prime}(\rho\sigma) [\sigma_{\rho}(\rho\sigma) + {\rm{i}}\tau_{\rho\theta}(\rho\sigma)] {\rm{d}}(\rho\sigma) \right\}} {{\rm{d}}\sigma}
\end{equation}
where
\begin{equation}
  \label{eq:5.3}
  \begin{aligned}
    \int z^{\prime}(\rho\sigma) [\sigma_{\rho}(\rho\sigma) + {\rm{i}}\tau_{\rho\theta}(\rho\sigma)] {\rm{d}}(\rho\sigma) =
    & \; \sum\limits_{\substack{ {k=-\infty} \\ {k \neq 0}}}^{\infty} {\rm{i}}A_{k-1}\frac{\rho^{k}}{k}\sigma^{k} - \sum\limits_{\substack{{k=-\infty} \\ {k \neq 0}}}^{\infty} {\rm{i}}B_{k-1} \frac{\rho^{-k}}{k} \sigma^{k} + (A_{-1}-B_{-1}){\rm{Ln}}\sigma \\
    & \; - \frac{z(\rho\sigma)-\overline{z(\rho\sigma)}}{\overline{z^{\prime}(\rho\sigma)}} \sum\limits_{k=-\infty}^{\infty} {\rm{i}}A_{-k} \rho^{-k} \sigma^{k} + C
  \end{aligned}
\end{equation}
wherein $C$ denotes an arbitrary complex integral constant. Similar to Eq.~(\ref{eq:4.18'}), Eq.~(\ref{eq:5.3}) can be rewritten as
\begin{equation}
  \label{eq:5.3'}
  \tag{5.3'}
  \begin{aligned}
    \int z^{\prime}(\rho\sigma) [\sigma_{\rho}(\rho\sigma) + {\rm{i}}\tau_{\rho\theta}(\rho\sigma)] {\rm{d}}(\rho\sigma)
    = & \; \sum\limits_{\substack{ {k=-\infty} \\ {k \neq 0}}}^{\infty} {\rm{i}}A_{k-1}\frac{\rho^{k}}{k}\sigma^{k} - \sum\limits_{\substack{{k=-\infty} \\ {k \neq 0}}}^{\infty} {\rm{i}}B_{k-1} \frac{\rho^{-k}}{k} \sigma^{k} + (A_{-1}-B_{-1}){\rm{Ln}}\sigma \\
    & \; + \sum\limits_{k=-\infty}^{\infty} \sum\limits_{l=-\infty}^{\infty} H_{l} \cdot {\rm{i}}A_{l-k} \rho^{l-k} \sigma^{k} + C \\
  \end{aligned}
\end{equation}
where
\begin{equation}
  \label{eq:5.4}
  \sum\limits_{k=-\infty}^{\infty} H_{k}\sigma^{k} \sim - \frac{z(\rho\sigma)-\overline{z(\rho\sigma)}}{\overline{z^{\prime}(\rho\sigma)}}
\end{equation}
The computation procedure of Eq.~(\ref{eq:5.4}) is completely the same to those in Appendix~\ref{sec:coefficients} for a given polar radius $\rho$, and would not be repeated here.

\begin{subequations}
  \label{eq:5.5}
  Eq.~(\ref{eq:5.1a}) can be transformed as
  \begin{equation}
    \label{eq:5.5a}
    \begin{aligned}
      \sigma_{\theta}(\rho\sigma) + \sigma_{\rho}(\rho\sigma)
      = 4\Re \left[ \frac{{\rm{i}}}{z^{\prime}(\rho\sigma) } \sum\limits_{k=-\infty}^{\infty} A_{k}\rho^{k}\sigma^{k} \right]
    \end{aligned}
  \end{equation}
  With Eq.~(\ref{eq:5.3'}), Eq.~(\ref{eq:5.1b}) can be further expressed as
  \begin{equation}
    \label{eq:5.5b}
    \sigma_{\rho}(\rho\sigma) + {\rm{i}}\tau_{\rho\theta}(\rho\sigma) = \frac{{\rm{i}}}{z^{\prime}(\rho\sigma)} \sum\limits_{k=-\infty}^{\infty} \left[ A_{k}\rho^{k} - B_{k}\rho^{-k-2} + (k+1) \sum\limits_{l=-\infty}^{\infty} H_{l}A_{l-k-1}\rho^{l-k-2} \right]\sigma^{k}
  \end{equation}
  Similarly, Eq.~(\ref{eq:5.1c}) can be simplified with Eq.~(\ref{eq:5.4}) as
  \begin{equation}
    \label{eq:5.5c}
    \begin{aligned}
      g(\rho\sigma) = 2G[u(\rho\sigma)+{\rm i}v(\rho\sigma)] =
      & \; {\rm{i}}\sum\limits_{k=1}^{\infty}\left( \kappa A_{k-1}\rho^{k} + B_{k-1}\rho^{-k} \right)\frac{\sigma^{k}}{k} - {\rm{i}}\sum\limits_{k=1}^{\infty}\left( \kappa A_{-k-1}\rho^{-k} + B_{-k-1}\rho^{k} \right)\frac{\sigma^{-k}}{k} + {\rm{i}}C_{0} \\
      & \; - {\rm{i}} \sum\limits_{k=-\infty}^{\infty} \sum\limits_{l=-\infty}^{\infty} H_{l} A_{l-k}\rho^{l-k}\sigma^{k} + {\rm{i}}(\kappa A_{-1} - B_{-1})\ln{\rho}
    \end{aligned}
  \end{equation} 
\end{subequations}
Eq.~(\ref{eq:5.5}) is the solution of the stress and displacement components in the mapping plane $ {\bm{\omega}} $ with an undetermined constant $ C_{0} $.

When $ \rho \rightarrow 1 $, Eq.~(\ref{eq:5.5}) would result in the stress and displacement along the outer boundary $ {\bm{c}}_{1} $. Note $ \overline{z(\sigma)} = \overline{z}(\sigma^{-1}) = z(\sigma) $ and $ \overline{z^{\prime}(\sigma)} = \overline{z^{\prime}}(\sigma^{-1}) = -\sigma^{2}z^{\prime}(\sigma) $ for $ \rho \rightarrow 1 $ according to Eq.~(\ref{eq:2.8}), Eq.~(\ref{eq:5.4}) would be zero, and subsequently Eq.~(\ref{eq:5.5}) can be reduced to the following expression as
\begin{subequations}
  \label{eq:5.6}
  \begin{equation}
    \label{eq:5.6a}
    \sigma_{\theta}(\sigma) + \sigma_{\rho}(\sigma) = 4\Re \left[ \frac{{\rm{i}}}{z^{\prime}(\sigma) } \sum\limits_{k=-\infty}^{\infty} A_{k}\sigma^{k} \right]
  \end{equation}
  \begin{equation}
    \label{eq:5.6b}
    \sigma_{\rho}(\sigma) + {\rm i}\tau_{\rho\theta}(\sigma) = \frac{{\rm{i}}}{z^{\prime}(\sigma)} \sum\limits_{k=-\infty}^{\infty} \left( A_{k} - B_{k} \right)\sigma^{k}
  \end{equation}
  \begin{equation}
    \label{eq:5.6c}
    g(\sigma) = 2G[u(\sigma)+{\rm i}v(\sigma)] = {\rm{i}}\sum\limits_{k=1}^{\infty}\left( \kappa A_{k-1} + B_{k-1}\right)\frac{\sigma^{k}}{k} - {\rm{i}}\sum\limits_{k=1}^{\infty}\left( \kappa A_{-k-1} + B_{-k-1} \right)\frac{\sigma^{-k}}{k} + {\rm{i}}C_{0}
  \end{equation} 
\end{subequations}
As mentioned in Eq.~(\ref{eq:4.11}), $ C_{0} $ can be determined by $ g(1) = 0 $ as
\begin{equation}
  \label{eq:5.7}
  C_{0} = \sum\limits_{k=1}^{\infty} \frac{1}{k} \left( \kappa A_{-k-1} + B_{-k-1} \right) - \sum\limits_{k=1}^{\infty} \frac{1}{k} \left( \kappa A_{k-1} + B_{k-1} \right)
\end{equation}
Eqs.~(\ref{eq:5.5c}) and~(\ref{eq:5.6c}) together indicate that displacement within geomaterial region in the mapping $ {\bm{\omega}} $ is always finite, and the original displacement singularity at infinity in the physical plane ($z\rightarrow\infty$) in Ref \cite[]{Zengguisen2019}, which is corresponding to point $\zeta = 1$ in the mapping plane, is eliminated.

Due to the sudden change of boundary condition along the outer boundary segments $ {\bm{c}}_{11} $ and $ {\bm{c}}_{12} $, the Gibbs phenomena would occur to cause stress and displacement oscillations, expecially along boundary $ {\bm{c}}_{1} $. The Lanczos filtering function \cite[]{lin2023reasonable,Lanczos1956} is applied to replace $ \sigma^{k} $ in the stress and displacement components with $ L_{k}\cdot \sigma^{k} $ to reduce such oscillations:
\begin{equation}
  \label{eq:5.8}
  L_{k} = \left\{
    \begin{aligned}
      & \; 1, \quad k = 0 \\
      & \; \sin\left( \frac{k}{N}\pi \right)/\left( \frac{k}{N}\pi \right), \quad {\rm{otherwise}}
    \end{aligned}
  \right.
\end{equation}
where $N$ takes the same value as that in Eqs.~(\ref{eq:4.4a'}) and~(\ref{eq:4.4b'}), indicating the same trunction is applied to Eqs.~(\ref{eq:5.5}) and~(\ref{eq:5.6}). It should be emphasized that the oscillations can only be reduced, but can not be fully eliminated.

The stress and displacement components caused by excavation in the physical plane can be obtained as
\begin{subequations}
  \label{eq:5.9} 
  \begin{equation}
    \label{eq:5.9a}
    \left\{
      \begin{aligned}
        & \sigma_{y}(z) + \sigma_{x}(z) = \sigma_{\theta}(\zeta) + \sigma_{\rho}(\zeta) \\
        & \sigma_{y}(z) - \sigma_{x}(z) + 2{\rm i}\tau_{xy}(z) = \left[ \sigma_{\theta}(\zeta) - \sigma_{\rho}(\zeta) + 2{\rm i} \tau_{\rho\theta}(\zeta) \right] \cdot \frac {\overline{\zeta}} {\zeta} \frac {\overline{z^{\prime}(\zeta)}} {z^{\prime}(\zeta)} \\
      \end{aligned}
    \right.
  \end{equation}
  \begin{equation}
    \label{eq:5.9b}
    u(z) + {\rm i}v(z) = u(\zeta) + {\rm i}v(\zeta)
  \end{equation}
\end{subequations}
where $ \sigma_{x}(z) $, $ \sigma_{y}(z) $, and $ \tau_{xy}(z) $ denote the horizontal, vertical, and shear stress components due to excavation in the lower half plane, respectively; $ u(z) $ and $ v(z) $ denote horizontal and vertical displacement components in the lower half plane, respectively. The final stress field within geomaterial in the rectangular coordinate system $ xOy $ in the physical plane is the sum of Eqs.~(\ref{eq:2.1}) and~(\ref{eq:5.9a}), and can be expressed as
\begin{equation}
  \label{eq:5.10}
  \left\{
    \begin{aligned}
      & \sigma_{x}^{\ast} = \sigma_{x} + \sigma_{x}^{0} \\
      & \sigma_{y}^{\ast} = \sigma_{y} + \sigma_{y}^{0} \\
      & \tau_{xy}^{\ast} = \tau_{xy} + \tau_{xy}^{0} \\
    \end{aligned}
  \right.
\end{equation}
Eq.~(\ref{eq:5.10}) can be transformed into polar form as
\begin{equation}
  \label{eq:5.11}
  \left\{
    \begin{aligned}
      & \sigma_{\theta}^{\ast}(\zeta) + \sigma_{\rho}^{\ast}(\zeta) = \sigma_{y}^{\ast}(z) + \sigma_{x}^{\ast}(z) \\
      & \sigma_{\theta}^{\ast}(\zeta) - \sigma_{\rho}^{\ast}(\zeta) + 2{\rm i} \tau_{\rho\theta}^{\ast}(\zeta) = \left[ \sigma_{y}^{\ast}(z) - \sigma_{x}^{\ast}(z) + 2{\rm i}\tau_{xy}^{\ast}(z) \right] \cdot \frac {\zeta} {\overline{\zeta}} \frac {z^{\prime}(\zeta)} {\overline{z^{\prime}(\zeta)}} \\
    \end{aligned}
  \right.
\end{equation}
where $\sigma_{\theta}^{\ast}(\zeta)$, $\sigma_{\rho}^{\ast}(\zeta)$, and $\tau_{\rho\theta}^{\ast}(\zeta)$ denote the final hoop, radial, and shear stress components in the polar coordinate system $ \rho o \theta $ in the mapping plane, respectively. Till now, the stress and displacement distributions in the geomaterial region $ {\bm \varOmega} $ are obtained, and our problem is solved.

\section{Solution verification}
\label{sec:verif}

In this section, the proposed solution in this paper would be verified via several numerical cases. The coefficients of the conformal mappings of these numerical cases are listed in Table~\ref{tab:1}. Cases 1 and 2 are shallow vertical-wall semicircle tunnels with different depths, and Cases 3 and 4 are shallow horseshoe tunnels with different depths. The tunnel shapes for all four cases above can be found in Fig.~\ref{fig:6}. Case 5 is a shallow circular tunnel with 5m radius and 10m depth. In all subsequent numerical cases, the geomaterial is plane-strain with $ \kappa = 3-4\nu $.

For numerical accuracy, the truncation number $N$ in Eqs~(\ref{eq:4.4a'}) and (\ref{eq:4.4b'}) and $M$ in Eqs.~(\ref{eqb:1}) and (\ref{eqb:2}) take 80 and 360, respectively. All the numerical cases are performed via the programming code {\tt{fortran}} of compiler {\tt{gcc-13.2.1}}, and the linear systems of iterations are solved using the {\tt{dgesv}} package of {\tt{lapack-3.11.0}}. All figures are drawn using {\tt{gnuplot-5.4}}. The {\tt{GNU General Public License}} is applied to all the codes within this paper, and the codes are released in author {\sf{Luobin Lin}}'s github repository {\emph{github.com/luobinlin987/noncircular-shallow-tunnelling-reasonable-displacement}}.

\begin{table}[htb]
  \centering
  \caption{Coefficients of the conformal mappings of the five cases for verification}
  \begin{threeparttable}
    \begin{tabular}{cccccc}
      \toprule
      Case & 1 & 2 & 3 & 4 & 5 \\
      \midrule
      $a$ & 3.628290 & 20.129710 & 3.5782933 & 20.038580 & 8.6602543 \\
      $r$ & 2.7780510E-1 & 5.4057535E-2 & 2.7684460E-1 & 5.3478912E-2 & 2.6794916E-1 \\
      $ b_{1} $ & 1.4741707E-2 & 2.4127561E-4 & -1.4576164E-2 & -6.1574816E-3 & - \\
      $ b_{2} $ & -3.7750916E-3 & -4.9415696E-4 & -8.8445228E-3 & -4.1469518E-4 & - \\
      $ b_{3} $ & -2.4711951E-3 & -3.0263293E-5 & -2.0210147E-3 & -4.7331423E-6 & - \\
      $ b_{4} $ & -7.8826493E-4 & -8.8660255E-7 & -2.2367365E-4 & 3.0791176E-7 & - \\
      $ b_{5} $ & -1.7980317E-4 & -5.7930598E-9 & 1.1456264E-5 & 9.9983977E-9 & - \\
      $ b_{6} $ & -3.0605208E-5 & 6.1401095E-10 & 6.3018406E-6 & -4.8124310E-11 & - \\
      $ b_{7} $ & -3.4798795E-6 & 2.8342474E-11 & 2.3385462E-6 & - & - \\
      $ b_{8} $ & -9.9314818E-8 & - & 3.4944773E-7 & - & - \\
      \bottomrule
    \end{tabular}
    \begin{tablenotes}
      \footnotesize
    \item Note: Cases 1-4 are cited from Ref \cite[]{Zengguisen2019}.
    \end{tablenotes}
  \end{threeparttable}
  \label{tab:1}
\end{table}

\subsection{Lanczos filtering and boundary conditions}
\label{sec:lanczos-filtering}

\begin{figure}[htb]
  \centering
  \includegraphics{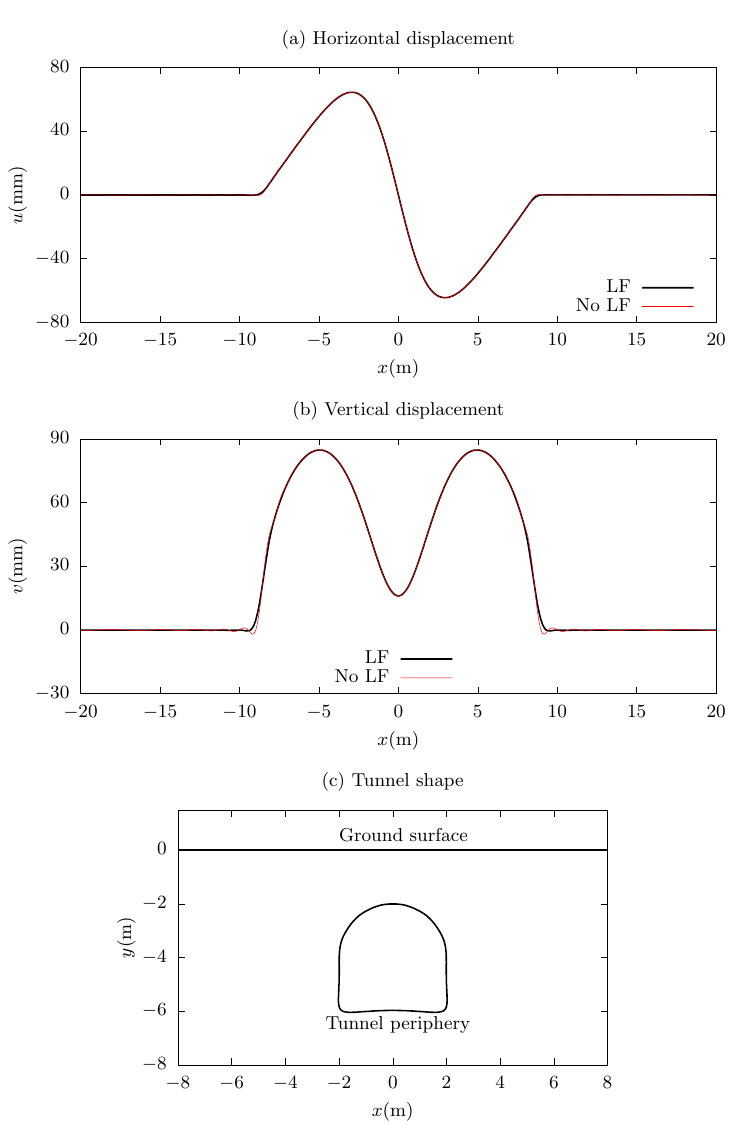}
  \caption{Displacement comparisons regarding of Lanczos filtering (Case 1, $\theta_{0} = 45^{\circ}, x_{0} \approx 8.75{\rm{m}}, \gamma = 20{\rm{kN/m^{3}}}, k_{0} = 0.8, \nu = 0.3, E = 1{\rm{MPa}}, N = 80$)}
  \label{fig:3}
\end{figure}

\begin{figure}[htb]
  \centering
  \includegraphics{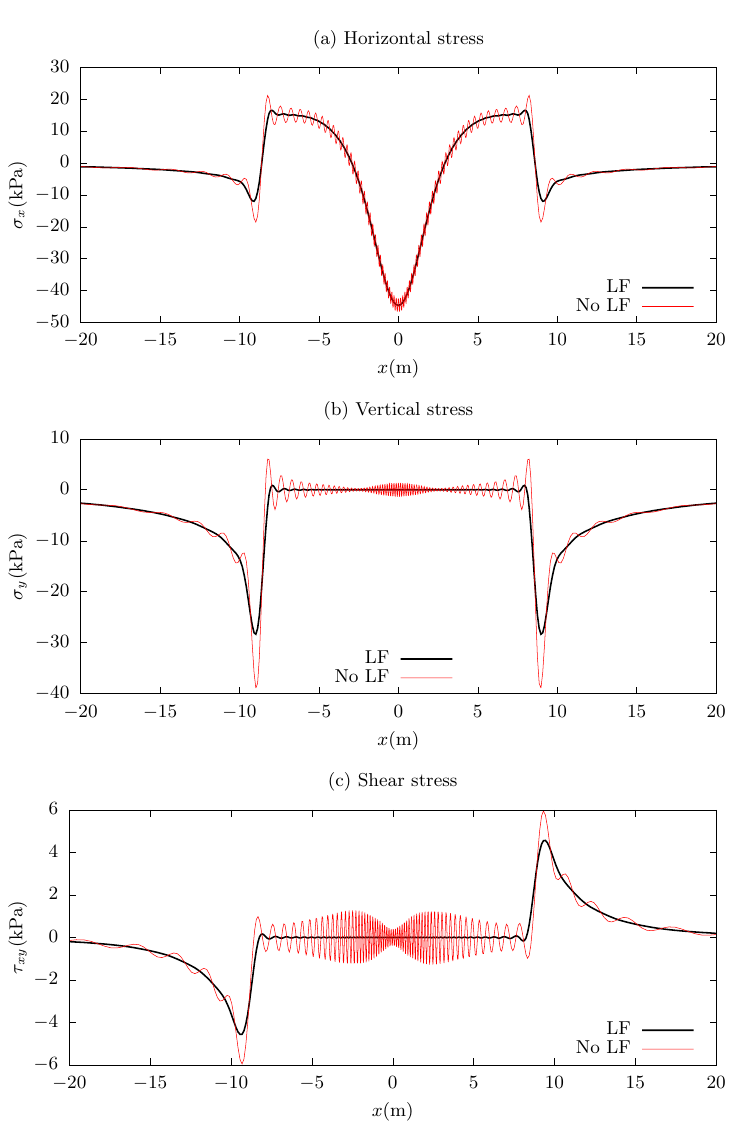}
  \caption{Stress comparisons regarding of Lanczos filtering (Case 1, $\theta_{0} = 45^{\circ}, x_{0} \approx 8.75{\rm{m}}, \gamma = 20{\rm{kN/m^{3}}}, k_{0} = 0.8, \nu = 0.3, E = 1{\rm{MPa}}, N = 80, q_{\max} = 43$)}
  \label{fig:4}
\end{figure}

The Lanczos filtering in Eq.~(\ref{eq:5.8}) is used to reduce the result oscillations caused by the Gibbs phenomena. To illustrate the neccessity of the Lanczos filtering, Case 1 is used to compare the results of stress and displacement components before and after using the Lanczos filtering. The other neccessary parameters are shown in the captions of Figs.~\ref{fig:3} and \ref{fig:4}. Substituting all the parameters and the coefficients of the conformal mappings of Case 1 into the proposed solution, the results are shown in Figs.~\ref{fig:3} and~\ref{fig:4}. The tunnel shape is shown in Fig.~(\ref{fig:3})c. The legend caption `LF` and `No LF` denote results with and without Lanczos filtering, respectively.

The red lines in Figs.~\ref{fig:3} and \ref{fig:4} show that both displacement and stress components along ground surface are oscillational without Lanczos filtering, which is identical to the foresight in Section~{\ref{sec:stress-displ-geom}}. The black lines, however, are more stable and with much less oscillation. Furthermore, the datum of the black lines are identical to the boundary conditions along the ground surface. To be specific, the following boundary conditions along ground surface in the rectangular coordinate system $xOy$ should be satisfied according to Eqs.~(\ref{eq:2.7a}) and (\ref{eq:2.7b}):
\begin{equation}
  \label{eq:6.1}
  \left\{
    \begin{aligned}
      & u = 0, v = 0, \quad |x| \geq x_{0} \\
      & \sigma_{y}^{\ast} = 0, \tau_{xy}^{\ast} = 0, \quad |x| \leq x_{0} \\
    \end{aligned}
  \right.
  , \quad
  {\rm{along \; ground \; surface}}
\end{equation}
Figs.~\ref{fig:3}a and~\ref{fig:3}b show that the rectangular displacement components $u$ and $v$ of the surface range $ x \geq |x_{0}| $ are zero, and Figs.~\ref{fig:4}b and~\ref{fig:4}c show that the surface stress components $ \sigma_{y} $ and $ \sigma_{xy} $ of the surface range $ x \leq |x_{0}| $ are zero. The results of these four subfigures are identical to the required boundary conditions along ground surface.

The reason that only the stress and displacement components along the ground surface are illustrated for comparisons before and after using the Lanczos filtering is that the Gibbs phenomena would be greatly and spontaneously reduced for regions away from the ground surface, even though the physical distance is very short. For convenience of practical computation, it can be very roughly interpreted that the Gibbs phenomena only occur along the ground surface.

The remaining boundary condition along tunnel periphery is also presented here according to Eq.~(\ref{eq:5.11}):
\begin{equation}
  \label{eq:6.2}
  \sigma_{\rho}^{\ast} = 0, \tau_{\rho\theta}^{\ast} = 0, \quad {\rm{along \; tunnel \; periphery}}
\end{equation}
Eq.~(\ref{eq:6.2}) is apparently reasonable, since the surface traction along tunnel periphery should be zero. To validate Eq.~(\ref{eq:6.2}), the following normalized stress components are computed in the subsequent numerical cases:
\begin{equation}
  \label{eq:6.3}
  \left\{
    \begin{aligned}
      & \Sigma_{\rho} = \frac{\max|\sigma_{\rho}^{\ast}|}{\gamma H} \\
      & \Sigma_{\rho\theta} = \frac{\max|\tau_{\rho\theta}^{\ast}|}{\gamma H}
    \end{aligned}
  \right.
\end{equation}
where $H$ denotes the depth of the tunnel.

\subsection{Constraining arc on solution convergence}
\label{sec:select-constr-arc}

\begin{figure}[htb]
  \centering
  \includegraphics{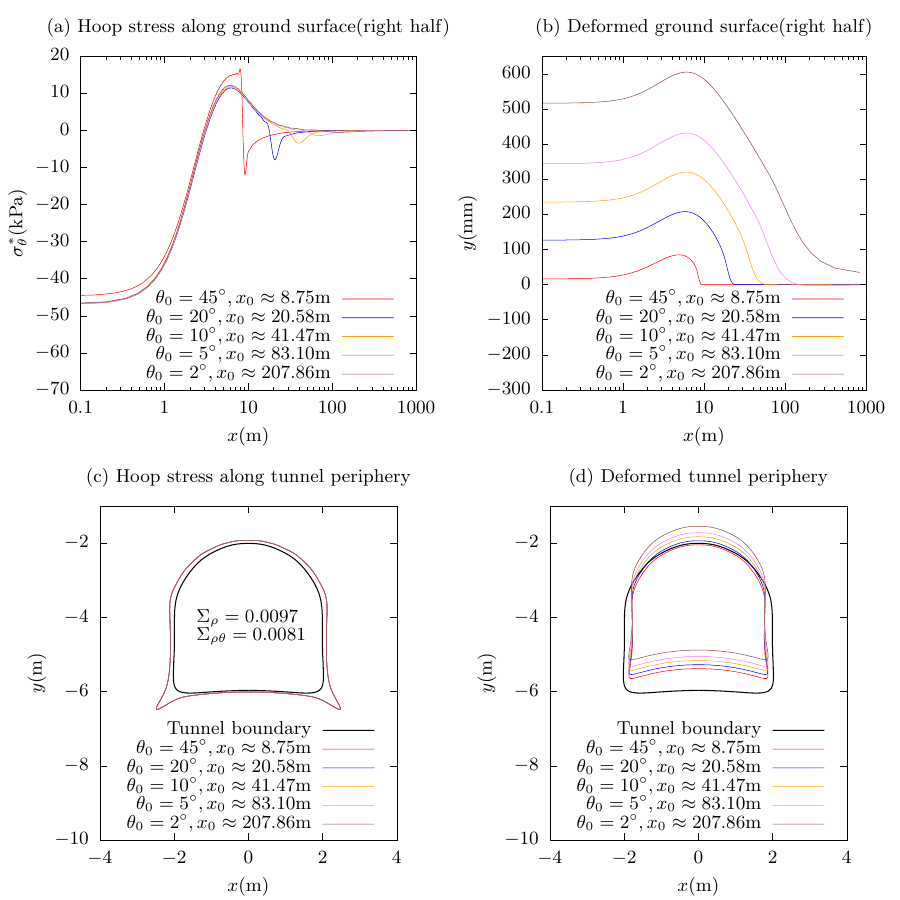}
  \caption{Stress and deformation comparisons on variation of constraining arc (Case 1, $\gamma = 20{\rm{kN/m^{3}}}, k_{0} = 0.8, \nu = 0.3, E = 1{\rm{MPa}}, N = 80$)}
  \label{fig:5}
\end{figure}

The range $ [-\theta_{0}, \theta_{0}] $ of the constraining arc $ \wideparen{t_{1}t_{2}} $ along the outer periphery in Fig.~\ref{fig:2}b would affect the solution convergence. The numerical cases in this subsection take the geometry of Case 1, and the other parameters are the same to the numerical case in last subsection, except for the constraining arc range $ \theta_{0} $. For comparisons, $ \theta_{0} $ takes the values of $ 45^{\circ} $, $ 20^{\circ} $, $ 10^{\circ} $, $ 5^{\circ} $, and $ 2^{\circ} $, respectively. Substituting all the neccessary parameters into the proposed solution, the stress and deformation along ground surface and tunnel periphery can be illustrated in Fig.~\ref{fig:5}. The horizontal axes of Figs.~\ref{fig:5}a and~\ref{fig:5}b apply logarithmic coordinate to cover larger range to better illustrate stress and deformation along ground surface.

Fig.~\ref{fig:5}a shows that as $ \theta_{0} $ gets smaller, the hoop stress along ground surface gradually become convergent, which seems to indicate that a smaller value of $ \theta_{0} $ would result in a more convergent solution. In theory, when $ \theta_{0} \rightarrow 0 $, the assumption of the ground surface of the proposed solution would be identical to the ones in Refs~\cite[]{Zengguisen2019,Self2020JEM,Lu2016,Verruijt_Strack2008buoyancy,Strack_Verruijt2002buoyancy,Strack2002phdthesis}. However, such assumptions are defective. The real-world geomaterial consists of soil particles, gravels, pebbles, and so on, thus, the deformation energy caused by excavation should be gradually consumed by all the nearby geomaterial constituents, and the far-field ground would remain undeformed. Therefore, it is more reasonable to keep the ground surface far from the tunnel to be constrained to absorb the possible overflowing energy which should have been consumed.

Fig.~\ref{fig:5}b indicates that a smaller value of $ \theta_{0} $ (or a larger range of free ground surface) would cause a larger upheaval of the ground surface, and all the upheavals share a similar trend of the surface deformation above the tunnel. Combining the hoop stress distribution in Fig.~\ref{fig:5}a, the value of $ \theta_{0} = 20^{\circ} $ would ensure both stress and displacement convergence for ground surface near the tunnel. As $ \theta_{0} $ gets larger, the displacement along ground surface would become inaccurate. To be specific, $ \theta_{0} = 10^{\circ} $, $ \theta_{0} = 5^{\circ} $, and $ \theta_{0} = 2^{\circ} $ should result in $ x_{0} = 41.47{\rm{m}} $, $ x_{0} = 83.10{\rm{m}} $, and $ x_{0} = 207.86{\rm{m}} $, respectively, indicating the ground surface outside of the coordinate $ (x_{0}, 0) $ should remain undeformed. However, it can be seen from Fig.~\ref{fig:5}b that the horizontal coordinates of the deformation curves for $ \theta_{0} = 10^{\circ} $, $ \theta_{0} = 5^{\circ} $, and $ \theta_{0} = 2^{\circ} $ exceed 50m, 100m, and 1000m, respectively, and are all larger than the expected values. Thus, Figs.~\ref{fig:5}a and~\ref{fig:5}b together indicate that the value $ \theta_{0} = 20^{\circ} $ would guarantee both stress and displacement convergence, as well as neccessary displacement boundary conditions along ground surface. The reason that a small value of $ \theta_{0} $, $ \theta_{0} = 2^{\circ} $ for instance, would cause overlarge displacement is a numerical one that the computation related to high-order $ {\rm{e}}^{\rm{i}k\theta_{0}} $ in the Taylor expansion in Eq.~(\ref{eq:4.3}) may contain significant errors.

Fig.~\ref{fig:5}c shows that hoop stress distribution and $ \Sigma_{\rho} $ and $ \Sigma_{\rho\theta} $ along tunnel periphery all remain stable for different $ \theta_{0} $, indicating that the boundary condition along tunnel periphery in Eq.~(\ref{eq:6.2}) is satisfied with small deviations and that the range of the constraining arc $ \wideparen{t_{1}t_{2}} $ would not greatly alter the stress distribution along tunnel periphery. Such a result is both analytically and numerically reasonable, since the stress distribution along tunnel periphery is the first derivative of Eq.~(\ref{eq:4.18}), whose right-hand side contains no coefficients related to $ \theta_{0} $. Thus, the variation of $ \theta_{0} $ should not affect the stress distribution along tunnel periphery.

Fig.~\ref{fig:5}d shows that the deformed tunnels with different $ \theta_{0} $ share a similar deformation pattern with different upheavals, which are in a similar pattern to Fig.~\ref{fig:5}b. Both Figs.~\ref{fig:5}b and~\ref{fig:5}d suggest that the variation of $ \theta_{0} $ would not greatly alter the relative deformation near the tunnel, and the existence of the constraining arc would constrain the far-field displacement of the geomaterial to be zero. However, the displacement solution in Ref~\cite[]{Zengguisen2019} would be infinitely large for far-field geomaterial, due to the logarithmic item related to the unbalanced resultant, just as analyzed in Section~\ref{sec:problem-2}. Moreover, the displacement solution in Ref~\cite[]{Zengguisen2019} need a `translation process` to eliminate the rigid-body displacement. Thus, the proposed solution eliminates the displacement singularity at infinity and the rigid-body displacement simultaneously, and is therefore more stable than the existing solution.

The above analyses suggest that the constraining arc range parameter $ \theta_{0} = 20^{\circ} $ can be used in the following numerical computation.

\subsection{Comparisons with Zeng's solution}
\label{sec:comp-with-zengs}

Though the problem presentation (balanced or unbalanced), the ground surface boundary condition, the to-be-determined variables, the solution procedure, and the final stress and displacement expressions between the proposed solution and Zeng's solution~\cite[]{Zengguisen2019} are quite different, both solutions deal with noncircular shallow tunnelling under initial stress field. To be specific, Eq.~(\ref{eq:4.15}) in this paper is the same to the right-hand side of Eq.~(28) in Zeng's solution~\cite[]{Zengguisen2019}, indicating that the boundary conditions along tunnel periphery between the proposed solution and Zeng's solution~\cite[]{Zengguisen2019} are the same. Owing to the uniqueness theorem of complex variable, the stress components along tunnel periphery between these two solution should be consequently the same for the same set of parameters.

To validate the proposed solution, four numerical cases (Cases 1-4 in Table~\ref{tab:1}) are computed for comparisons with Zeng's solution~\cite[]{Zengguisen2019}. To be identical, the parameters take the same values as Zeng's solution by setting $ \gamma = 27{\rm{kN/m^{3}}} $ and $ \nu = 0.3 $. The geometries of Cases 1 and 2 in Table~\ref{tab:1} are respectively the same to those two comparison cases of Fig.5 in Zeng's solution~\cite[]{Zengguisen2019}, and the lateral coefficient takes $ k_{0} = 0.25 $ as well. Meanwhile, the geometries of Cases 3 and 4 are respectivley the same to the numerical examples of Figs.7a and 7d in Zeng's solution~\cite[]{Zengguisen2019}, and the lateral coefficient takes $ k_{0} = 1.5 $ for the largest hoop stress cases. Under such parameter selections in the above four cases, the stress components along each tunnel periphery of all Cases 1-4 solved by the proposed solution should be the same to the corresponding ones already solved in Zeng's solution~\cite[]{Zengguisen2019}.

\begin{figure}[htb]
  \centering
  \includegraphics{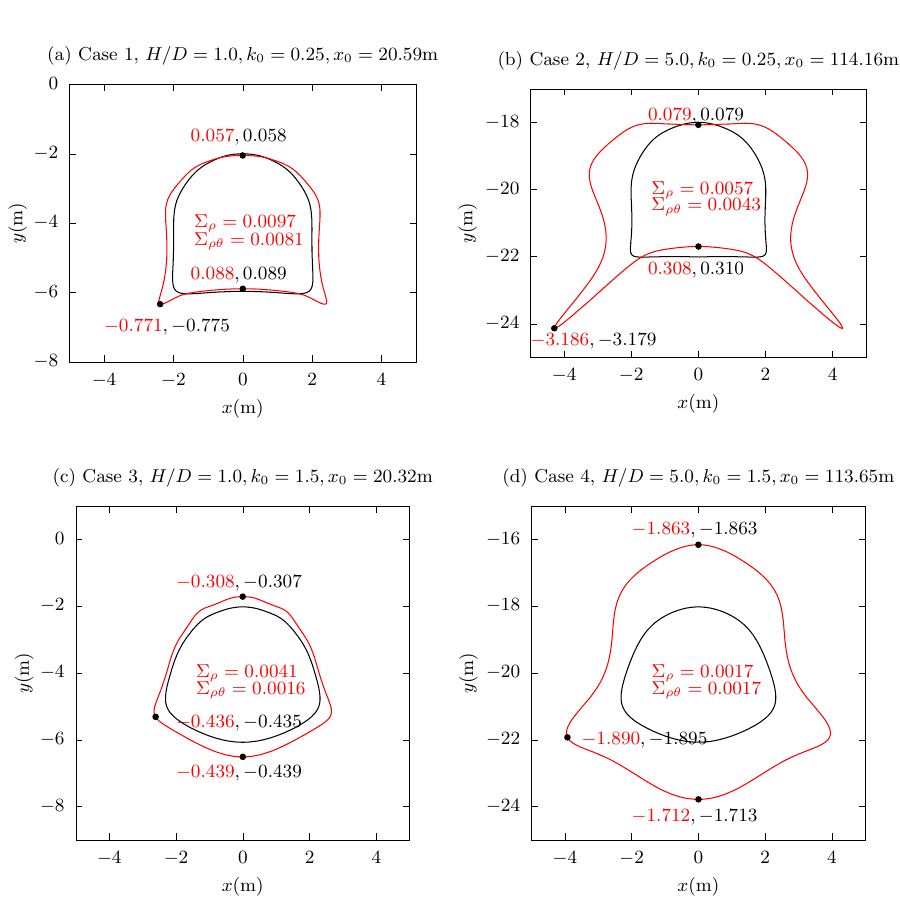}
  \caption{Comparisons of hoop stress with Zeng's solution \cite[]{Zengguisen2019}($\theta_{0} = 20^{\circ}, \gamma = 27{\rm{kN/m^{3}}}, \nu = 0.3, N = 80$)}
  \label{fig:6}
\end{figure}

Repectively substituting the above parameter sets into Cases 1-4 gives the stress component comparisons in Fig.~\ref{fig:6}. In Fig.~\ref{fig:6}, the black and red lines denote tunnel peripheries and hoop stress distributions for Cases 1-4, respectively, and the black dots denote the hoop stress values of several specific points along tunnel peripheries of Cases 1-4, which can be found in the corresponding figures in Zeng's solution~\cite[]{Zengguisen2019}. The red and black numbers are the computed results of the hoop stresses by the proposed solution and the existing ones cited from the corresponding figures in Zeng's solution~\cite[]{Zengguisen2019}, respectively. The red lines of hoop stress distributions and the black dots of hoop stress comparisons of specific points along tunnel periphery together indicate that the hoop stresses between these two solutions are very approximate. The computed results of $ \Sigma_{\rho} $ and $ \Sigma_{\rho\theta} $ are less than $ 10^{-2} $, and further indicate that the traction-free tunnel periphery is satisified, as mentioned in Eq.~(\ref{eq:6.2}). The stress comparisons in Fig.~\ref{fig:6} validate the proposed solution for arbitrary tunnel shapes.

\subsection{Comparisons with Lin's solution}
\label{sec:comp-with-lins}

Since last section indicates that the proposed solution can be applicable to arbitrary tunnel shapes, the proposed solution should be applicable to circular tunnel as well. By simply setting all $ b_{k} = 0 $, the conformal mapping in Eq.~(\ref{eq:2.8}) is degenerated into the Verrujt's conformal mapping\cite[]{Verruijt1997traction,Verruijt1997displacement}. Correpondingly, the stress and displacement solution can be degenerated, and should be identical to the one in Ref \cite[]{lin2023reasonable}. In this subsection, we would conduct a numerical case to compare the degenerated version of the proposed solution with Lin's solution \cite[]{lin2023reasonable} for further verification.

\begin{figure}
  \centering
  \includegraphics{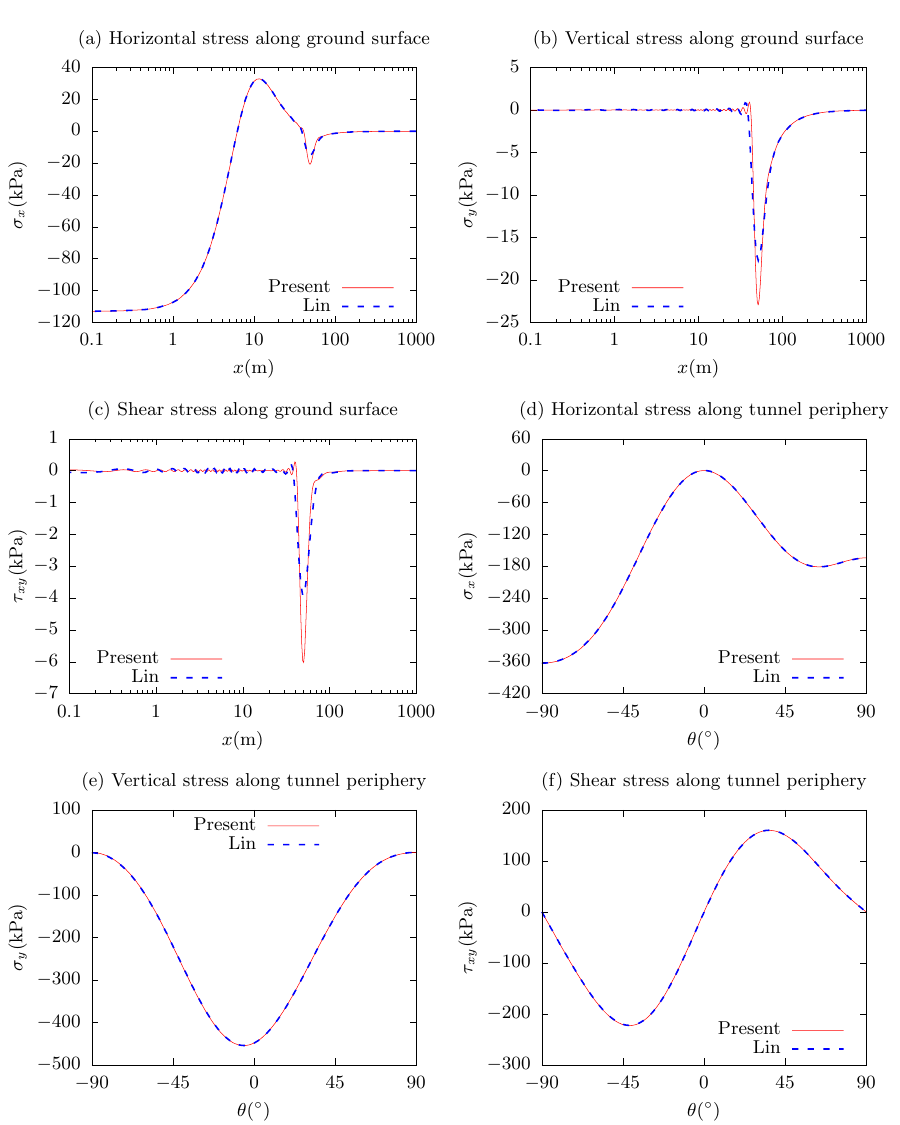}
  \caption{Comparisons of stress components with Lin's solution \cite[]{lin2023reasonable}(Case 5, $\theta_{0} = 20^{\circ}, x_{0} \approx 49.11{\rm{m}}, \gamma = 20{\rm{kN/m^{3}}}, k_{0} = 0.8, \nu = 0.3, E = 20{\rm{MPa}}, N = 80$)}
  \label{fig:7}
\end{figure}

\begin{figure}[htb]
  \centering
  \includegraphics{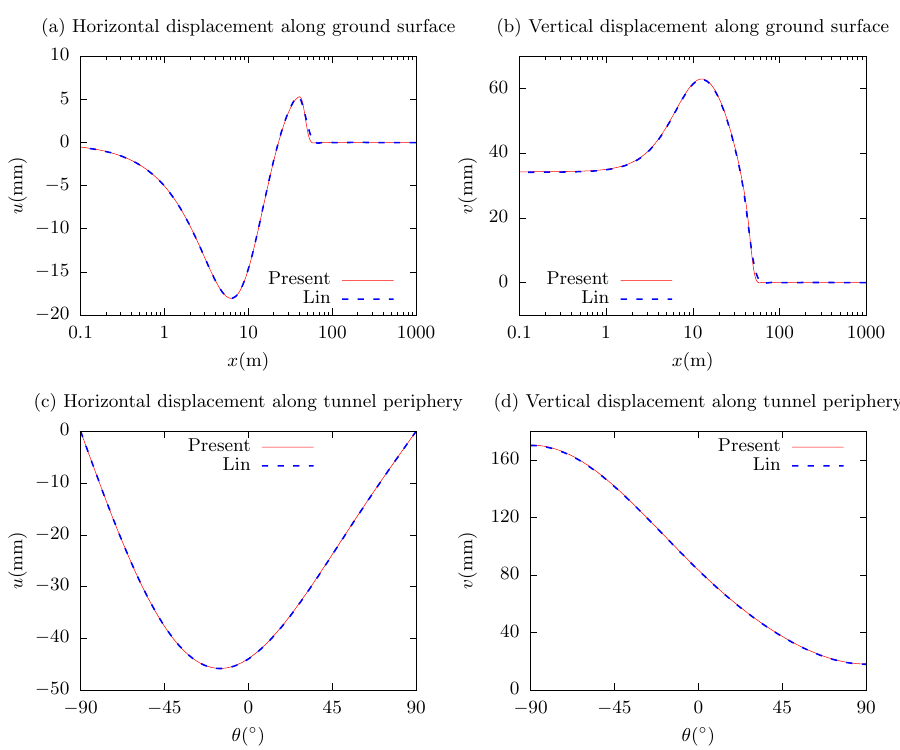}
  \caption{Comparisons of displacement components with Lin's solution \cite[]{lin2023reasonable}(Case 5, $\theta_{0} = 20^{\circ}, x_{0} \approx 49.11{\rm{m}}, \gamma = 20{\rm{kN/m^{3}}}, k_{0} = 0.8, \nu = 0.3, E = 20{\rm{MPa}}, N = 80$)}
  \label{fig:8}
\end{figure}

The parameters used in this numerical case are listed in the caption of Fig.~\ref{fig:7}, and the tunnel geometry takes the conformal mapping coefficients of Case 5 in Table~\ref{tab:1}. The programming codes of Lin's solution \cite[]{lin2023reasonable} can be found in the corresponding github repository along with that paper. Substituting all the parameters into these two solutions gives Figs.~\ref{fig:7} and~\ref{fig:8}. The datum in Figs.~\ref{fig:7} and~\ref{fig:8} all suggest highly identity of both stress and displacement components for these two solutions. Besides, the normalized radial and shear stress along tunnel periphery of the proposed solution are $ \Sigma_{\rho} = 0.0012 $, $ \Sigma_{\rho\theta} = 0.0004 $, which are also less than $ 10^{-2} $ and indicate that the traction-free tunnel periphery is satisfied.

\section{Stress and deformation}
\label{sec:stress-deformation}

Last section detailedly verifies the proposed solution and its convergence, and focuses more on the mathematical aspect. In this section, more general numerical cases focusing on tunnel engineering are conducted to investigate the stress and deformation distribution within geomaterial caused by arbitrary tunnel excavation considering initial stress field. The tunnel shapes are square, horseshoe, vertical-wall semicircle, and vertical horseshoe, and the corresponding coefficients of conformal mapping that control shallow tunnel shape are the same to those in Table 1 of Ref~\cite[]{Zengguisen2019}. The unit weight of geomaterial is chosen as $ \gamma = 20{\rm{kN/m^{3}}} $. The elastisity modulus and Poisson's ratio are chosen as $ E = 20{\rm{MPa}} $ and $ \nu = 0.3 $, and the plane strain condition is considered with $ \kappa = 3-4\nu $. The constraining arc takes $ \theta_{0} = 20^{\circ} $, as verified and recommended in last section. The truncation numbers take $ M = 360 $ and $ N = 80 $. The lateral coefficient is set to be a variable with values of $ k_{0} = 0.4, 0.8, 1.2, 1.6 $ to investigate its effects on stress and deformation distribution.

\subsection{Stress distribution}
\label{sec:stress-distribution}

\begin{figure}[htb]
  \centering
  \includegraphics{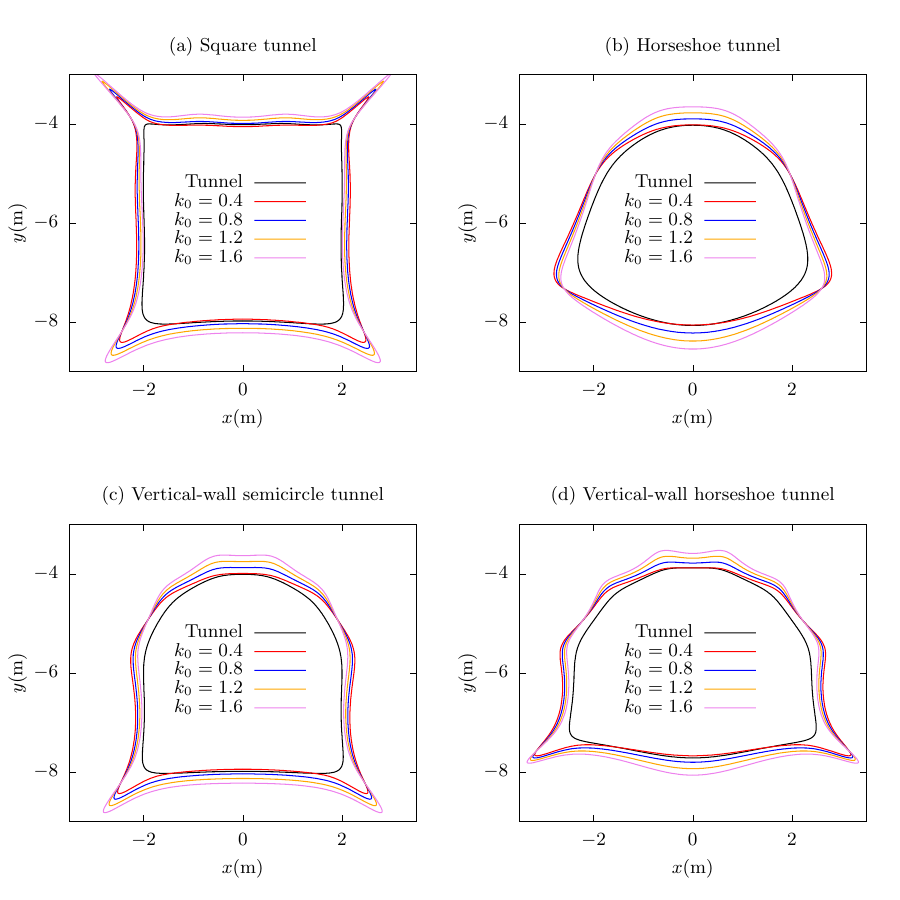}
  \caption{Hoop stress distribution along tunnel periphery for different tunnel shapes subjected to different lateral coefficients ($ \theta_{0} = 20^{\circ} $, $ \gamma = 20{\rm{kN/m^{3}}} $, $ E = 20{\rm{MPa}} $, $ \nu = 0.3 $, $ N = 80 $)}
  \label{fig:9}
\end{figure}

Substituting all the neccessary parameters mentioned above into the proposed solution, the hoop stress distributions along tunnel periphery for all four tunnel shapes and four lateral coefficients are obtained in Fig.~\ref{fig:9}. Most visually and intuitively, Fig.~\ref{fig:9} shows that serious stress concentrations exist for the corners along tunnel periphery, especially for the square tunnel. These stress concentrations indicate that tunnel corners should be avoided or smoothened if possible to geometrically reduce the abrupt change of stress along tunnel periphery for tunnel safety, and to more reasonably use the mechanical capacity of the nearby geomaterial.

Fig.~\ref{fig:9} also shows that when the lateral coefficient gets larger, the hoop stresses around tunnel top, bottom, and corners increase remarkably, while the ones around tunnel side walls get smaller with a much less numerical decrease. Such a stress variation around different tunnel segements indicates that the lateral coefficient would primarily affect the top and bottom of a shallow tunnel. Thus, possible engineering methods may be neccessary to protect tunnel top and bottom for stratum of high lateral coefficient. A more detailed inspection suggests that tensile stresses possibly exist around tunnel top and bottom for all four tunnel shapes when lateral coefficient is small $ k_{0} = 0.4 $. Correpondingly, tunnel top and bottom may also be reinforced for stratum of low lateral coefficient.

\subsection{Deformation distribution}
\label{sec:deform-distr}

\begin{figure}[htb]
  \centering
  \includegraphics{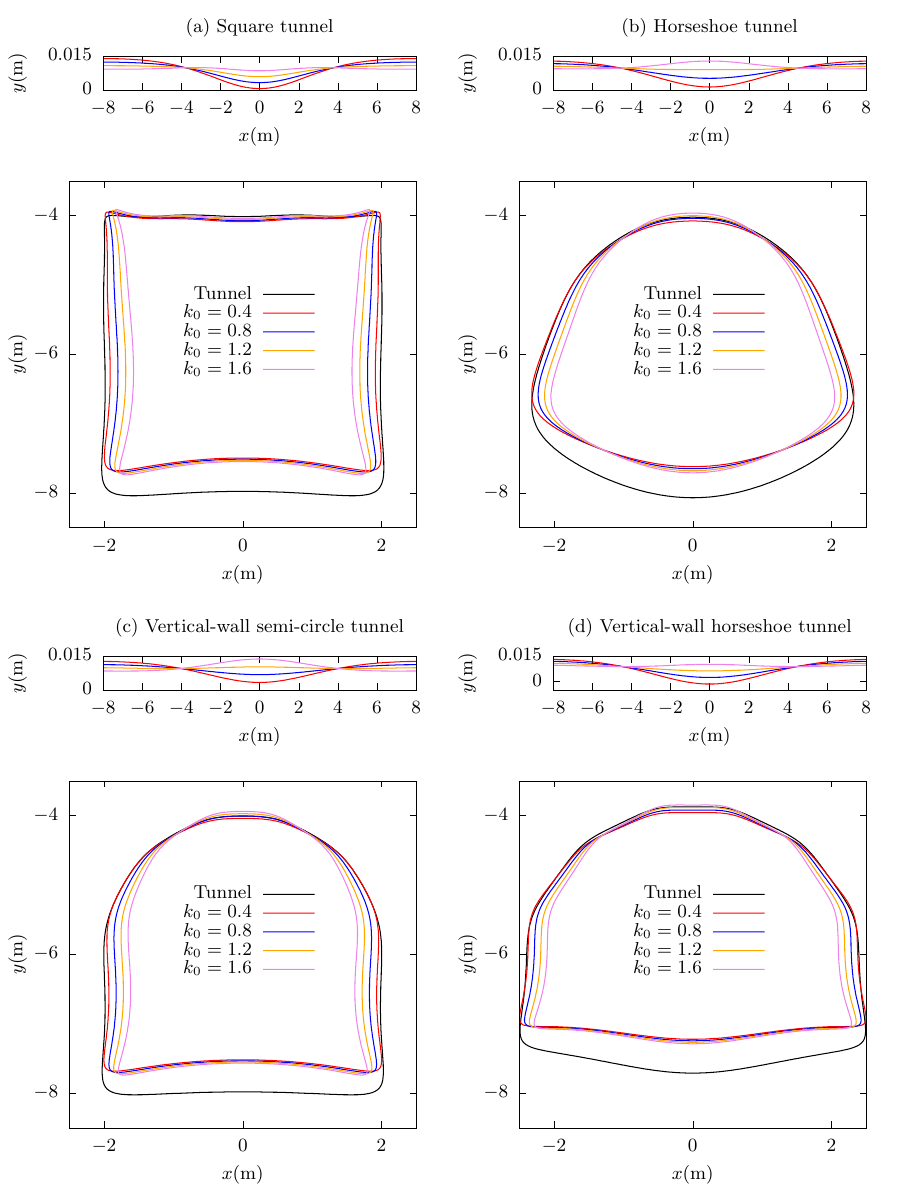}
  \caption{Deformation along ground surface (as is) and tunnel periphery (10 times magnified for better illustration) for different tunnel shapes subjected to different lateral coefficients ($ \theta_{0} = 20^{\circ} $, $ \gamma = 20{\rm{kN/m^{3}}} $, $ E = 20{\rm{MPa}} $, $ \nu = 0.3 $, $ N = 80 $)}
  \label{fig:10}
\end{figure}

The deformations along ground surface and tunnel periphery for all four tunnel shapes and four lateral coefficients are shown in Fig.~\ref{fig:10}. For better illustration, the deformations along tunnel periphery are magnified by 10 times, while the ones along ground surface remain as is.

Fig.~\ref{fig:10} shows that as the lateral coefficient increases from 0.4 to 1.6, the side walls of all four tunnels gradually shrink inwards, while the ground surfaces gradually turn from settlement troughs to corresponding upheavals. The reason is that when $ k_{0} $ is larger than 1, the horizontal stress would be larger than the vertical one, and the geomaterial would spontaneously and horizontally squeezes towards the exavated cavity. Therefore, side wall deformation should be monitored during tunnel excavation to avoid unfavorable engineering hazards.

Fig.~\ref{fig:10} also shows that all four tunnels share a generally upward lift or buoyancy. The reason is that the geomaterial excavation provides an upward resultant along tunnel periphery, as shown in Eq.~(\ref{eq:2.3}), and the geomaterial around the tunnel would be consequently lifted upwards. The solution in this paper assumes that geomaterial is elastic, and the excavation is considered instantaneous. However, tunnel excavation is  generally time-dependent, and the stress equilbriums during excavation may be achieved multiple times, and the upward lift in real tunnel engineering should be generally less than the results obtained via the proposed solution. However, in other words, the deformation solved via the proposed solution may be treated as a possible limit case for the worst situation.

\section{Conclusions}
\label{sec:conclusion}

The proposed solution deals with stress and displacement of a new mechanical model aiming at eliminating the far-field displacement singularity by applying analytic contituation to tranform the mixed boundaries to a homogenerous Riemann-Hilbert problem, which is solved via an iterative procedure. The subsequent verifications all suggest good agreements with expected boundary conditions and existing analytical solutions. Several numerical cases are conducted to investigate the stress and deformation distribution along the ground surface and tunnel periphery and the following conclusions are drawn:

(1) Serious stress concentrations occurring around corners along tunnel periphery should be avoided by smoothening tunnel geometry to reduce possible damage for tunnel safety.

(2) Shallow tunnel top and bottom may be more vulnerable than side walls to lateral coefficient variation of complicated stratum encountered in tunnel engineering.

(3) Deformation monitoring may be neccessary for shallow tunnelling in stratum of relatively high lateral coefficient, due to noteworthy inward deformation along tunnel side walls.

(4) A general upward lift may occur along shallow tunnel periphery and nearby geomaterial due to excavation, which should be considered and computed in design and before construction, and possible engineering measures may be planned beforehand to deal with such an unfavorable buoyancy.

\section{Further discussions}
\label{sec:futher-discussion}

Though the proposed solution successfully eliminates the far-field displacement singularity to simultaneously obtain correct stress and reasonable deformation for noncircular shallow tunnelling, some defects still exist, and should be improved in further studies:

(1) The conformal mapping in Eq.~(\ref{eq:2.8}) is only a backward one, which maps the mapping plane onto the physical plane, while the corresponding forward conformal mapping, which maps the physical plane onto the mapping plane, is still not found. Such a lack of forward conformal mapping makes us have to take the polar coordinates in the mapping plane as variables. Thus, we can not directly input an arbitrary rectangular coordinate $ (x,y) $ in the physical plane into the proposed solution to compute the corresponding stress and deformation.

(2) The conformal mapping in Eq.~(\ref{eq:2.8}) is only suitable for axisymmetrical noncircular cavity, and can not deal with asymmetrical cavity.

(3) The solution process of the coefficients of the conformal mapping in Eq.~(\ref{eq:2.8}) is dependent on solving a non-convex objective function~\cite[]{Zengguisen2019}, and a heavy parallel computation may be needed, which requires parallel computation infrastructures (high-performance GPUs for instance) and more sophisticated coding skills ({\tt{CUDA}} for instance). One possible parallel solution for the conformal mapping using particle swarm method~\cite[]{kennedy1995particle} and {\tt{openacc}} coding technique can be found in Ref~\cite[]{lin2022solution}. Both hardware and skill requirements elevate the threshold to obtain the correct conformal mapping.

Therefore, an interesting and promising study direction is to find a new pair of forward and backward conformal mappings that mutually maps a lower half plane with an asymmetrical cavity and a unit annulus onto each other. The theoretical basis to construct such a pair of conformal mappings should be more solid; the solution method to obtain such a pair of conformal mappings should better be a linear system instead of non-convex optimization; and the computation of the solution method should better be very fast.

\clearpage
\section*{Acknowledgement}
\label{sec:acknowledgement}

This study is financially supported by the Natural Science Foundation of Fujian Province, China (Grant No. 2022J05190), the Scientific Research Foundation of Fujian University of Technology (Grant No. GY-Z20094 and GY-H-22050), and the National Natural Science Foundation of China (Grant No. 52178318). The authors would like to thank Professor Changjie Zheng, Ph.D. Yiqun Huang, and Associate Professor Xiaoyi Zhang for their suggestions on this study.

\clearpage
\appendix

\section{Expansion of Eq.~(\ref{eq:4.15'})}
\label{sec:expansion}

\begin{equation}
  \label{eqa:1}
  - \gamma \int_{A}^{B} y{\rm d}x = {\rm{i}} \gamma \sum\limits_{j=1}^{3} \int_{A}^{B} C_{j}(\sigma) {\rm d}\sigma
\end{equation}
where
\begin{subequations}
  \label{eqa:2}
  \begin{equation}
    \label{eqa:2a}
    C_{1}(\sigma) =  - a^{2}r \left[ \frac{1}{(1-r\sigma)^{3}} + \frac{r}{(\sigma-r)^{3}} + \frac{1-r^{2}}{(1-r\sigma)^{2}(\sigma-r)^{2}} \right]
  \end{equation}
  \begin{equation}
    \label{eqa:2b}
    \begin{aligned}
      C_{2}(\sigma) = & \; \frac{a}{4} \sum\limits_{k=1}^{K} kb_{k}(r^{k}+r^{-k}) \left( \frac{1+r\sigma}{1-r\sigma}\sigma^{k-1} + \frac{1+r\sigma}{1-r\sigma}\sigma^{-k-1} + \frac{\sigma+r}{\sigma-r}\sigma^{k-1} + \frac{\sigma+r}{\sigma-r}\sigma^{-k-1} \right) \\
      + & \; \frac{a}{2} \sum\limits_{k=1}^{K} b_{k}(r^{k}-r^{-k}) \left[ \frac{r\sigma^{k}}{(1-r\sigma)^{2}} + \frac{r\sigma^{-k}}{(1-r\sigma)^{2}} + \frac{r\sigma^{k}}{(\sigma-r)^{2}} + \frac{r\sigma^{-k}}{(\sigma-r)^{2}} \right]
    \end{aligned}
  \end{equation}
  \begin{equation}
    \label{eqa:2c}
    C_{3}(\sigma) = -\frac{1}{4} \sum\limits_{k=1}^{K}\sum\limits_{l=1}^{K} b_{k}(r^{k}-r^{-k}) \cdot lb_{l}(r^{l}+r^{-l}) \cdot (\sigma^{k+l-1} + \sigma^{-k-l-1} + \sigma^{k-l-1} + \sigma^{-k+l-1})
  \end{equation}
\end{subequations}
Therefore, the indefinite integrals can be correspondingly obtained as
\begin{equation}
  \label{eqa:3}
  \int_{A}^{B} C_{j}(\sigma) {\rm d}\sigma = E_{j}(\sigma) + D_{j} \cdot {\rm Ln}\sigma, \quad j = 1,2,3
\end{equation}
where
\begin{subequations}
  \label{eqa:4}
  \begin{equation}
    \label{eqa:4a}
    E_{1}(\sigma) = - \frac{a^{2}}{2(1-r\sigma)^{2}} + \frac{a^{2} r^{2}}{2(\sigma-r)^{2}} + \frac{a^{2} r}{(1-r^{2}) (\sigma-r)} - \frac{a^{2} r^{2}}{(1-r^{2}) (1-r\sigma)} - \frac{2a^{2}r^{2}}{(1-r^{2})^{2}}\ln\frac{1-r\sigma^{-1}}{1-r\sigma}
  \end{equation}
  \begin{equation}
    \label{eqa:4b}
    \begin{aligned}
      E_{2}(\sigma)
      = & \; a\sum\limits_{k=1}^{K} b_{k}(r^{2k}-r^{-2k})\left( \frac{1}{1-r\sigma} - \frac{r}{\sigma-r} \right) + a\sum\limits_{k=1}^{K}b_{k} \sum\limits_{l=1}^{k} \frac{k-l}{l}(r^{2k-l}-r^{-2k+l})(\sigma^{l}-\sigma^{-l}) \\
      + & \; a\sum\limits_{k=1}^{K} kb_{k}(r^{2k}+r^{-2k})\ln\frac{1-r\sigma^{-1}}{1-r\sigma} + \frac{a}{2} \sum\limits_{k=1}^{K} b_{k}(r^{k}+r^{-k}) \left( \frac{1}{1-r\sigma} + \frac{r}{\sigma-r} \right)(\sigma^{k}-\sigma^{-k})
    \end{aligned}
  \end{equation}
  \begin{equation}
    \label{eqa:4c}
    \begin{aligned}
      E_{3}(\sigma) = & \; -\frac{1}{4} \sum\limits_{k=1}^{K} b_{k}(r^{k}-r^{-k}) \sum\limits_{ \substack{ l=1 \\ l \neq k} }^{K} lb_{l}(r^{l}+r^{-l}) \cdot \left( \frac{\sigma^{k+l}}{k+l} - \frac{\sigma^{-k-l}}{k+l} + \frac{\sigma^{k-l}}{k-l} + \frac{\sigma^{-k+l}}{-k+l} \right) \\
      & \; - \frac{1}{8}\sum\limits_{k=1}^{K} b_{k}^{2}(r^{2k}-r^{-2k})(\sigma^{2k}-\sigma^{-2k})
    \end{aligned}
  \end{equation}
\end{subequations}
\begin{subequations}
  \label{eqa:5}
  \begin{equation}
    \label{eqa:5a}
    D_{1} = - \frac{2a^{2}r^{2}}{(1-r^{2})^{2}}
  \end{equation}
  \begin{equation}
    \label{eqa:5b}
    D_{2} = 2a \sum\limits_{k=1}^{K} kb_{k}r^{2k}
  \end{equation}
  \begin{equation}
    \label{eqa:5c}
    D_{3} = -\frac{1}{2} \sum\limits_{k=1}^{K} kb_{k}^{2} (r^{2k}-r^{-2k})
  \end{equation}
\end{subequations}
where the following integrations are applied:
\begin{subequations}
  \label{eqa:6}
  \begin{equation}
    \label{eqa:6a}
    \int \frac{1+r\sigma}{1-r\sigma}\sigma^{k-1}{\rm{d}}\sigma = -\frac{\sigma^{k}}{k} + \frac{2}{k}\frac{1}{1-r\sigma}\sigma^{k} - \frac{2r}{k} \int \frac{\sigma^{k}}{(1-r\sigma)^{2}}{\rm{d}}\sigma + 2r^{-1} \cdot \delta(k,1), \quad k \geq 1
  \end{equation}
  \begin{equation}
    \label{eqa:6b}
    \int \frac{1+r\sigma}{1-r\sigma}\sigma^{-k-1}{\rm{d}}\sigma = \frac{\sigma^{-k}}{k} - \frac{2}{k}\frac{1}{1-r\sigma}\sigma^{-k} + \frac{2r}{k} \int \frac{\sigma^{-k}}{(1-r\sigma)^{2}}{\rm{d}}\sigma, \quad k \geq 1
  \end{equation}
  \begin{equation}
    \label{eqa:6c}
    \int \frac{\sigma+r}{\sigma-r}\sigma^{k-1}{\rm{d}}\sigma = \frac{\sigma^{k}}{k} + \frac{2}{k}\frac{r}{\sigma-r}\sigma^{k} + \frac{2r}{k} \int \frac{\sigma^{k}}{(\sigma-r)^{2}}{\rm{d}}\sigma - 2r \cdot \delta(k,1), \quad k \geq 1
  \end{equation}
  \begin{equation}
    \label{eqa:6d}
    \int \frac{\sigma+r}{\sigma-r}\sigma^{-k-1}{\rm{d}}\sigma = -\frac{\sigma^{-k}}{k} - \frac{2}{k}\frac{r}{\sigma-r}\sigma^{-k} - \frac{2r}{k} \int \frac{\sigma^{-k}}{(\sigma-r)^{2}}{\rm{d}}\sigma, \quad k \geq 1
  \end{equation}
\end{subequations}
where
\begin{subequations}
  \label{eqa:7}
  \begin{equation}
    \label{eqa:7a}
    \int \frac{\sigma^{k}}{(1-r\sigma)^{2}}{\rm{d}}\sigma = \sum\limits_{l=1}^{k} \frac{k-l}{l} r^{-k-1+l}\sigma^{l} + \frac{r^{-k-1}}{1-r\sigma} + kr^{-k-1}\ln(1-r\sigma), \quad k \geq 1
  \end{equation}
  \begin{equation}
    \label{eqa:7b}
    \int \frac{\sigma^{-k}}{(1-r\sigma)^{2}}{\rm{d}}\sigma = -\sum\limits_{l=1}^{k-1} \frac{l}{k-l} r^{l-1}\sigma^{-k+l} + \frac{r^{k-1}}{1-r\sigma} - kr^{k-1}\ln(1-r\sigma) + kr^{k-1}{\rm{Ln}}\sigma, \quad k \geq 1
  \end{equation}
  \begin{equation}
    \label{eqa:7c}
    \int \frac{\sigma^{k}}{(\sigma-r)^{2}}{\rm{d}}\sigma = \sum\limits_{l=1}^{k} \frac{k-l}{l} r^{k-l-1}\sigma^{l} - \frac{r^{k}}{\sigma-r} + kr^{k-1}\ln(1-r\sigma^{-1}) + kr^{k-1}{\rm{Ln}}\sigma, \quad k \geq 1
  \end{equation}
  \begin{equation}
    \label{eqa:7d}
    \int \frac{\sigma^{-k}}{(\sigma-r)^{2}}{\rm{d}}\sigma = -\sum\limits_{l=1}^{k-1} \frac{l}{k-l} r^{-l-1} \sigma^{-k+l} - \frac{r^{-k}}{\sigma-r} - kr^{-k-1}\ln(1-r\sigma^{-1}), \quad k \geq 1
  \end{equation}
\end{subequations}
\begin{equation*}
  \delta(k,1) = \left\{
    \begin{aligned}
      & 1, \quad k = 1 \\
      & 0, \quad {\rm{otherwise}} \\
    \end{aligned}
  \right.
\end{equation*}
Eq.~(\ref{eqa:7}) comes from Ref~\cite[]{Zengguisen2019}. 

\clearpage

\section{Coefficients in Eq.~(\ref{eq:4.19})}
\label{sec:coefficients}

Since both expressions in front of the symbol $\sim$ are analytic and smooth, the following matchimg point method can be used to approximately solve the coefficients in Eq.~(\ref{eq:4.19}). We can uniformly select $2M+1$ points as $ \theta_{i} = \frac{i-1}{2M+1}2\pi $, and then $ \sigma_{i} = {\rm e}^{{\rm i}\theta_{i}} $. Eq.~(\ref{eq:4.19}) ensure the following equations stand:
\begin{subequations}
  \label{eqb:1}
  \begin{equation}
    \label{eqb:1a}
    \sum\limits_{k=-M}^{M} I_{k} \sigma_{i}^{k} = - \frac{z(r\sigma_{i})-\overline{z(r\sigma_{i})}}{\overline{z^{\prime}(r\sigma_{i})}}, \quad i \in [1, 2M+1] \cap {\bm Z}
  \end{equation}
  \begin{equation}
    \label{eqb:1b}
    \sum\limits_{k=-M}^{M} J_{k} \sigma_{i}^{k} = F(\sigma_{i}) +  \sum\limits_{j=1}^{4} E_{j}(\sigma_{i}), \quad i \in [1, 2M+1] \cap {\bm Z}
  \end{equation}
\end{subequations}
The coefficients can be solved via two definite linear systems. To guarantee numerical stability of the iteration of Eq.~(\ref{eq:4.24}), $ I_{k} = 0 $ is mandatorily set without loss of accuracy, if $ |I_{k}| \leq \delta $, where the threshold $\delta$ can be selected as $ 10^{-12} $ for instance. Subsequently, the infinite bilateral series in Eq.~(\ref{eq:4.19}) can be simulated by the finite series in Eq.~(\ref{eqb:1}) as
\begin{subequations}
  \label{eqb:2}
  \begin{equation}
    \label{eqb:2a}
    \sum\limits_{k=-\infty}^{\infty} I_{k} \sigma^{k} \sim \sum\limits_{k=-M}^{M} I_{k} \sigma^{k}
  \end{equation}
  \begin{equation}
    \label{eqb:2b}
    \sum\limits_{k=-\infty}^{\infty} J_{k} \sigma^{k} \sim \sum\limits_{k=-M}^{M} J_{k} \sigma^{k}
  \end{equation}
\end{subequations}

\clearpage

\clearpage

\clearpage

\clearpage

\clearpage

\clearpage

\clearpage

\clearpage

\clearpage

\clearpage

\clearpage
\bibliographystyle{plainnat}
\bibliography{./manu7}  %%% Uncomment this line and comment out the ``thebibliography'' section below to use the external .bib file (using bibtex) .

\end{document}